\documentclass[a4paper,11pt,reqno]{amsart}
\usepackage[latin1]{inputenc}
\usepackage[english]{babel}
\usepackage{amssymb,amsfonts, amsmath, color}
\usepackage[margin=2.7cm]{geometry}
\usepackage{xcolor}
\usepackage{graphicx, color, enumerate}
\usepackage[latin1]{inputenc}
\usepackage[active]{srcltx}
\usepackage{tikz}
\usepackage{pgf}
\usepackage{etex}
\usepackage{verbatim}
\usepackage{tikz-3dplot}
\usepackage{amsmath}
\usepackage{amsfonts}
\usepackage{amssymb}
\usepackage{amsthm}
\usepackage{algpseudocode}
\usepackage{float}
\usepackage{xcolor}
\usepackage{mathrsfs}
\usepackage{import}
\usepackage{geometry}
\usepackage{fancyhdr}
\usepackage{fp}

\newtheorem{theorem}{Theorem}[section]
\newtheorem{proposition}[theorem]{Proposition}
\newtheorem{corollary}[theorem]{Corollary}
\newtheorem{lemma}[theorem]{Lemma}
\newtheorem{remark}[theorem]{Remark}
\newtheorem{definition}[theorem]{Definition}

\newcommand{\bcl}{\begin{center}}
\newcommand{\ecl}{\end{center}}
\newcommand{\brl}{\begin{right}}
\newcommand{\erl}{\end{right}}
\newcommand{\ben}{\begin{enumerate}}
\newcommand{\een}{\end{enumerate}}
\newcommand{\overliner}{\begin{array}}
\newcommand{\earr}{\end{array}}
\newcommand{\btab}{\begin{tabular}}
\newcommand{\etab}{\end{tabular}}
\newcommand{\bdoc}{\begin{document}}
\newcommand{\edoc}{\end{document}}
\newcommand{\beqy}{\begin{eqnarray}}
\newcommand{\eeqy}{\end{eqnarray}}

\newcommand{\beqi}{\begin{eqnarray*}}
\newcommand{\eeqi}{\end{eqnarray*}}
\newcommand{\bitem}{\begin{itemize}}

\newcommand{\eitem}{\end{itemize}}
\newcommand{\nln}{\newline}
\newcommand{\newt}{\newtheorem}

\newcommand{\pa}{\partial}
\newcommand{\re}{{I\!\!R}}
\newcommand{\Rn}{\R^N}
\newcommand{\xr}{x\in\R }
\newcommand{\x}{\times}
\newcommand{\dyle}{\displaystyle}
\newcommand{\ene}{{I\!\!N}}
\newcommand{\irn}{\int\limits_{\R^N}}
\newcommand{\io}{\int\limits_{\O}}
\newcommand{\meas}{{\rm meas\,}}
\newcommand{\dif}{\nabla_{xy}}
\newcommand{\sign}{{\rm sign}}
\newcommand{\map}{\longrightarrow }
\newcommand{\imp}{\Longrightarrow }
\renewcommand{\div}{\nabla\cdot }
\newcommand{\sen}{{\rm sen\,}}
\newcommand{\tg}{{\rm tg\,}}
\newcommand{\arcsen}{{\rm arcsen\,}}
\newcommand{\arctg}{{\rm arctg\,}}
\newcommand{\supp}{{\textsl supp\ }}
\newcommand{\ity}{\int_{-\iy}^{+\iy}}
\newcommand{\limit}{\lim\limits}
\newcommand{\limi}{\limit_{n\to\infty}}
\newcommand{\sumi}{\sum\limits_{n=1}^{\infty}}
\newcommand{\ulu}{\underline u}
\newcommand{\ulw}{\underline w}
\newcommand{\ulz}{\underline z}
\newcommand{\ulv}{\underline v}
\newcommand{\uls}{\underline s}
\newcommand{\olu}{\overline u}
\newcommand{\olv}{\overline v}
\newcommand{\ols}{\overline s}
\newcommand{\ob}{\overline\b}
\newcommand{\ovar}{\overline\var}
\newcommand{\wv}{\widetilde v}
\newcommand{\wu}{\widetilde u}
\newcommand{\ws}{\widetilde s}
\renewcommand{\a }{\alpha }
\renewcommand{\b }{\beta }
\newcommand{\g }{\gamma}
\newcommand{\G }{\mathcal G }
\renewcommand{\d }{\delta }

\newcommand{\D }{\Delta }
\newcommand{\e }{\varepsilon }
\newcommand{\z }{\zeta }
\renewcommand{\l }{\lambda }
\renewcommand{\L }{\Lambda }
\newcommand{\m }{\mu }
\newcommand{\n }{\nabla }
\newcommand{\s }{\sigma }
\newcommand{\Sig }{\Sigma }
\renewcommand{\t }{\tau }
\newcommand{\var }{\varphi }
\renewcommand{\o }{\omega }
\renewcommand{\O }{\Omega }
\newcommand{\R}{{\mathbb{R}}}
\newcommand{\N}{{\mathbb{N}}}
\newcommand{\bC}{{\bf C}}
\newcommand{\bZ}{{\bf Z}}
\newcommand{\bN}{{\bf N}}
\newcommand{\bQ}{{\bf Q}}
\newcommand{\bK}{{\bf K}}
\newcommand{\bI}{{\bf I}}
\newcommand{\bv}{{\bf v}}
\newcommand{\bV}{{\bf V}}
\DeclareMathOperator{\suppo}{supp} \DeclareMathOperator{\di}{div}

\newenvironment{Proof}{\Rmovelastskip\vskip12pt
plus 1pt \noindent\em\rm}{\hfill {\qed \hskip .2cm}}

\begin{document}

\title[]{Uniqueness of solutions to \\ elliptic and parabolic equations \\ on metric graphs}
\author{Giulia Meglioli}

\address{\hbox{\parbox{5.7in}{\medskip \noindent{Giulia Meglioli, \\Fakult\"at f\"ur Mathematik, \\Universit\"at Bielefeld, \\33501, Bielefeld, Germany \\ [3pt] \emph{E-mail address: }{\tt gmeglioli@math.uni-bielefeld.de}}}}}

\author{Fabio Punzo}

\address{\hbox{\parbox{5.7in}{\medskip \noindent{Fabio Punzo, \\Dipartimento di Matematica, \\Politecnico di Milano, \\Piazza Leonardo da Vinci 32, 20133, Milano, Italy \\ [3pt] \emph{E-mail address: }{\tt fabio.punzo@polimi.it}}}}}


\subjclass[2010]{}

\maketitle

\begin{abstract}
We investigate uniqueness of solutions to certain classes of elliptic and parabolic equations posed on {\it metric graphs}. 
In particular, we address the linear Schr\"odinger equation with a potential, and the heat equation with a variable density. We assume suitable growth conditions on the solutions, which are related to the behaviour at infinity of the potential or of the density.
\end{abstract}

\section{Introduction}\label{sec0}
We are concerned with both elliptic equations of the type
\begin{equation}\label{problema}
\Delta u -V u=0 \quad \text{ in }\; \G,
\end{equation}
and parabolic equations of the form
\begin{equation}\label{problema2}
\rho\partial_tu-\Delta u=0 \quad \text{ in }\; \G\times(0,T),
\end{equation}
completed with the initial condition
\begin{equation}\label{incond}
u(x,0)=0\quad\text{ for}\,\,\,x\in\G\,.
\end{equation}
Here $\mathcal G$ is a {\it metric} graph over an infinite weighted graph $(\mathcal V, E, l)$ with set of vertices and edges $\mathcal V$ and $E$, respectively, and weight function $l$; $\Delta$ stands for the Laplace operator on $\mathcal G$ coupled with Kirchoff boundary conditions on the vertices. Furthermore, $V:\mathcal G\to \mathbb (0, +\infty)$ and $\rho:\mathcal G\to (0, +\infty)$ are given functions, we refer to $V$ as a {\it potential} and to $\rho$ as a {\it density}.

The goal of this paper is to investigate uniqueness of solutions to \eqref{problema} and to \eqref{problema2}-\eqref{incond}, under suitable growth conditions.

\smallskip

The study of uniqueness for solutions to \eqref{problema} or to \eqref{problema2}-\eqref{incond}, posed on $\mathbb R^N$, but also on Riemannian manifolds, has a very long history in the literature, also including more general operators, see e.g. \cite{AB, BMeP, EKP, F, Gaf, Grig1, Gri, Grig, IKO, L, MePu, MR, PT, Punzo2, T}. Very recently, such uniqueness question has also been investigated on {\it discrete graphs}, see e.g. \cite{BMP, BP, Huang, HKS, Meg, MPgraph}. Discrete (or combinatorial) graphs are composed of sets of vertices connected by edges. These edges primarily serve as abstract relationships between vertices or as carriers of supplementary attributes such as weights and directions. Thus a function is said to be defined on a discrete graph whenever it is defined only on its vertices; therefore, also the gradient and the Laplace operator take into account only the values of the function on the vertices. In contrast, metric graphs are conceived as spatially continuous networks, where edges are treated as physical line segments joined at vertices. This continuity enables dynamic phenomena to evolve along the edges, making them suitable for modeling spatial systems found in physics and biology, such as river networks, cave formations, and vascular systems (see \cite{SCA, SMA, SFPP}). More broadly, graphs emerge in various scientific disciplines, including quantum physics, chemistry, engineering, mathematical biology, and mathematical fields such as algebra and number theory (see \cite{BF, BK, MSV, R}). Within this context, dynamic behaviors on graphs are often described through differential equations posed on the edges, subject to specific boundary conditions imposed at the vertices (see \cite{Mu2}).

In last years, also the investigation of elliptic and parabolic problems on {\it metric} graphs is attracting the attention of various authors, see \cite{AST, BDS, CJS1, DpR, DST, EFMM, KMN, PSV, PT, PuTe2, PuTe3, vb, vb1}.

\subsection{Outline of our results}
In this paper, we prove the uniqueness of solutions in suitable weighted $L^p$ spaces for both the elliptic equation \eqref{problema} and the Cauchy problem \eqref{problema2}-\eqref{incond}, in the setting of a \textit{metric graph}. Fix a reference point $x_0\in \mathcal V$. We denote by $d(x, x_0)$ the distance between any $x\in \mathcal G$  and $x_0$. About the potential $V$, we either assume that $V$ is bounded away from zero, that is
\begin{equation}\label{Vbound}
V(x)\ge V_0\quad\text{for all}\,\,\,x\in\G\,,
\end{equation}
for some $V_0>0$,
or that, $V$ decays to zero as $d(x, x_0)\to \infty$ with an appropriate rate, in the sense that
\begin{equation}\label{Vbound2}
V(x)\ge V_0\,[d(x,x_0)+k]^{-\theta}\quad\text{for all}\,\,\,x\in\G\,,
\end{equation}
for some $V_0>0$, $0<\theta\le2$ and $k\ge1$. Similarly, we either assume that
\begin{equation}\label{rhobound}
\rho(x)\ge \rho_0\quad\text{for all}\,\,\,x\in\G\,,
\end{equation}
for some $\rho_0>0$,
or that
\begin{equation}\label{rhobound2}
\rho(x)\ge \rho_0[d(x,x_0)+k]^{-\theta}\quad\text{for all}\,\,\,x\in\G\,.
\end{equation}
for some $\rho_0>0$, $0<\theta\le2$ and $k\ge1$.
More precisely, consider first equation \eqref{problema}. We differentiate the results based on the behavior of $V$.
\begin{itemize}
\item If $V$ is a positive function bounded from below by a positive constant $V_0$, as described in \eqref{Vbound}, then we show that the only solution to \eqref{problema} in the space $L^p_{\varphi_\beta}(\G)$, for all $p\ge1$, is $u\equiv 0$. Here,
$$
\varphi_\beta(x):=e^{-\beta d(x,x_0)}\quad \text{for any}\,\,\,x\in\G,
$$
for some $\beta>0$. The cases $p\ge2$ and $1\le p<2$ are treated separately, see Theorem \ref{teo1}. In particular, for $1\le p<2$ an additional assumption on the underlying graph, namely \eqref{degAss}, is required to obtain our result.
\item If $V$ is a positive function that may vanish at infinity (see \eqref{Vbound2}), then we show that $u\equiv 0$ is the only solution to \eqref{problema} in the space $L^p_{\varphi_\lambda}(\G)$ for some $p\ge1$ with
$$
\varphi_\lambda(x):=[d(x,x_0)+k]^{-\lambda}\quad \text{for any}\,\,\,x\in\G,
$$
where $\lambda>0$ and $k\ge1$. This is further detailed in Theorem \ref{teo3}. Again, the case $1\le p<2$ requires a more accurate analysis and the extra assumption \eqref{degAss} on the graph structure.
\end{itemize}
A similar analysis is conducted for the Cauchy problem \eqref{problema2}-\eqref{incond}. 
\begin{itemize}
\item If $\rho$ is a positive function bounded from below by a positive constant $\rho_0$, see \eqref{rhobound}, then we show that the only solution to \eqref{problema2}-\eqref{incond} in the space $L^p_{\varphi_\beta}(\G\times(0,T))$ for some $p\ge1$, is $u\equiv 0$. The weight $\varphi_\beta$ is the same function defined above. As for the elliptic equation, the cases $p\ge2$ and $1\le p<2$ need to be treated separately, see Theorem \ref{teo5}. In fact, for $1\le p<2$, the extra assumption \eqref{degAss} on the underlying graph must hold.
\item If $\rho$ is a positive function which may decay at infinity (see \eqref{rhobound2}), then we show that $u\equiv 0$ is the only solution to \eqref{problema2}-\eqref{incond} in $L^p_{\varphi_\lambda}(\G\times(0,T))$ for all $p\ge1$, as detailed in Theorem \ref{teo7}. As in the previous cases, for $1\le p<2$, a more accurate analysis is needed, in fact, this case requires the extra assumption \eqref{degAss} on the underlying graph.
\end{itemize}

\subsection{Organization of the paper}
The paper is organized as follows. In Section \ref{mf} we provide an overview on the metric graphs, and we give the definitions of solutions to both equation \eqref{problema} and problem \eqref{problema2}-\eqref{incond}. Section \ref{statement} contains all our main results divided between the elliptic and the parabolic cases. Afterwards, in Section \ref{sec2} we state and prove some auxiliary results. Finally the main theorem concerning the elliptic equation is proved in Section \ref{sec5} and the one concerning the parabolic problem is shown in Section \ref{sec6}.

\section{Mathematical framework}\label{mf}
In this section, we gather some basic notions and definitions related to metric graphs, along with some preliminaries concerning analysis on metric graphs.

\subsection{The metric graph setting} 
Given a countable set of {\it vertices}, a countable set of {\it edges}, and a function $l: E\to (0, +\infty)$, usually called {\it weight}, consider the {\it weighted} graph $(\mathcal V, E, l)$. Now, we regard the edges as segments glued together at the vertices, and $l(e)$ as the length of the edge $e$. Let
\[\mathcal E:=\prod_{e\in E} \{e\}\times (0, l(e)).\]
We may give the following definition
\begin{definition}
The metric graph $\mathcal G$ over the weighted graph $(\mathcal V, E, l)$ is the pair $(\mathcal V, \mathcal E)$. 
\end{definition}
We endow the metric graph $\mathcal G$ with the maps $i:E\to \mathcal V$ and $j:\{e\in E:\,l_e<\infty\}\to\mathcal  V$, which define the initial and the final point of an edge. Both the initial and final points of an edge are called endpoints of an edge. An oriented weighted graph uniquely determines a metric graph once we adopt the convention that $i(e)=0$ and $j(e)=l(e).$ 

In the sequel, by an abuse of notation, we will generally refer to points in $\G$ as $x\in\G$ when either $x=v\in \mathcal V$ or $x\in \{e\}\times (0, l(e))$ for some $e\in E$. Sometimes, for simplicity, we identify $e\equiv (0, l(e)).$ Thus, we also write $x\in e$ or $x\in (0, l(e)).$ Therefore, roughly speaking, the metric graph $\mathcal G$ can also be regarded as $\mathcal V\cup E$. The same consideration for the point $x$ holds when we consider a function $u:\mathcal G\to \mathbb R.$

\begin{definition}\label{degree}
We define the {\it degree} $\operatorname{deg}_v\in\N$ of a vertex $v\in \mathcal V$ as the number of edges $e\ni v$. The {\it inbound degree} $\operatorname{deg}_v^+$ (respectively, the {\it outbound degree} $\operatorname{deg}_v^-$) is the number of edges with $j(e)=v$ (respectively, with $i(e)=v$). Clearly, $\operatorname{deg}_v=\operatorname{deg}_v^++\operatorname{deg}_v^-$.
\end{definition}
\smallskip

Additionally, we give the following definitions:
\begin{itemize}
\item A metric graph $\mathcal G$ is {\it finite} if both $E$ and $\mathcal V$ are finite sets; conversely, we say that the graph is {\it infinite};
\item Let $v,w\in \mathcal V$. A path connecting $v$ and $w$ is a set $\{x_1,...,x_n\}\subseteq\mathcal G$, ($n\in\N$) such that $x_1=v$, $x_n=w$ and for all $k=1,...,n-1$ there exists an edge $e_k$ such that $x_k,x_{k+1}\in\overline{\mathcal G}_{e_k}$. If $v\equiv w$ then the path is closed. A closed path is called cycle if it does not pass through the same node more than once.
\item A metric graph $\mathcal G$ does not posses {\it loops} if no edge connects a vertex to itself, and no two vertices are connected by more than one edge, i.e. $i(e)\ne j(e)$ for all $e\in E$.
\item A metric graph $\mathcal G$ is {\it connected} if, for any two vertices $v,w\in \mathcal V$, there exists a path joining $v$ to $w$. A connected graph without cycles is a {\it tree}.
\item A metric graph $\mathcal G$ is {\it locally finite} if $\operatorname{deg}_v<\infty$ for all $v\in \mathcal V$.
\item The {\it boundary} of the metric graph is $\partial\mathcal G:=\{v\in\mathcal  V\,:\, \operatorname{deg}_v=1\}$.
\end{itemize}

It is easily seen that a connected metric graph $\mathcal G$ is a metric measure space. In fact, every two points $x,y\in\mathcal G$ can be regarded as vertices of a path $P=\{x\equiv x_1,x_2,...,x_n\equiv y\}$ connecting them. The length of the path is the sum of the $n$ edges $e_k$, i.e. $l(P):= \sum_{k=1}^n l(e_k)\,,$ and the distance $d=d(x,y)$ between $x$ and $y$ is
$$
d(x,y)\,:=\,\inf\,\{l(P)\,|\, \text{$P$ connects $x$ and $y$}\}\,.
$$
Therefore $\mathcal G$ is a metric space, with Borel $\s$-algebra $\mathcal{B}=\mathcal{B}(\G)$. A Radon measure $\m:\mathcal{B}\mapsto [0,\infty]$ is induced on $\mathcal G$ by the Lebesgue measure $\l$ on each interval $I_e$, namely
\begin{equation}\label{demu}
\m(\Omega)\,:=\, \sum_{e\in E} \l(I_e\cap \Omega) \quad\text{for any $\Omega\in\mathcal{B}$}\,.
\end{equation}

In the sequel we always make the following hypothesis:
\begin{equation}\label{ipotesigrafo}
\begin{cases}
&\mathcal G \text{ is an infinite, connected, locally finite metric graph without loops}, \\ 
& l_e<\infty \text{ for any } e\in E\,,\\ 
& \text{there exists a vertex } x_0\in \mathcal V \text{ and a sequence }
\{v_n\}\subset \mathcal V \text{ s.t. } d(x_0, v_n)\to +\infty\,.
\end{cases}
\end{equation}

\subsection{Functional spaces on metric graphs}\label{funcSpace}
Set $I_e:=(0, l(e)).$ We denote by $\mathfrak F$ the set of all functions $f:\mathcal G\to \mathbb R.$ Given any $f\in \mathfrak F$, let
\[f_e:=f|_{\bar I_e}.\]
Therefore every function $f\in\mathfrak F$ canonically induces a countable family of functions $\{f_e\}_{e\in E}$, $f_e:\bar I_e\to\R$, thus we write $$f=\bigoplus_{e\in E}f_e\,.$$

In view of \eqref{demu}, for any measurable $f\in\mathfrak F$ we set
\begin{equation}\label{int1}
\begin{aligned}
&\int_{\mathcal G} f\,d\m\,:=\, \sum_{e\in E} \int_{I_e} f_e\,dx\,, \\
& \int_{\Omega} f\,d\m\,:=\, \int_{\Omega} f\chi_\Omega\,dx
\quad\text{for any $\Omega\in\mathcal{B}(\G)$\,,}
\end{aligned}
\end{equation}
where $\chi_{\Omega}$ denotes the characteristic function of the set $\Omega$ and the usual notation $dx\equiv d\l$ is used. Accordingly, for any $p\in[1,\infty]$ the Lebesgue spaces $L^p(\G)\equiv L^p(\mathcal G,\m)$ are immediately defined:
$$
L^p(\mathcal G)\,:=\,\bigoplus_{e\in E}\,L^p((0,l(e)),\l)
$$
with norm
$$
\begin{aligned}
&\|f\|_p\,:=\,\sum_{e\in E}\|f_e\|_p\, =\, \sum_{e\in E} \left(\int_0^{l_e} |f_e|^p\,dx\right)^{\frac1p} \;\; \text{if $p\in [1,\infty)$}\,, \\
& \|f\|_\infty\,:=\,{\rm ess}\, \sup_{x\in\mathcal G}\, |f(x)|\,.
 \end{aligned}
$$
Moreover, let $\varphi:\mathcal G\to\R$, $\varphi>0$ be continuous. We define, for each $p\in[1,+\infty)$, the {\it weighted Lebesgue spaces} $L^p_\varphi(\G)\equiv L^p_\varphi(\mathcal G,\m)$ as follows
$$
L^p_\varphi(\mathcal G)\,:=\,\bigoplus_{e\in E}\,L^p_\varphi((0,l(e)),\l)=\left\{f\in\mathfrak F\,:\, \sum_{e\in E} \left(\int_0^{l_e} |f_e|^p\varphi\,dx\right)^{\frac1p}<+\infty\right\}.
$$
We say that $u$ is continuous in $\mathcal G$, and we write $u\in C(\G)$, whenever $u_e\in C(\bar I_e)$, for any $e\in E$, and at any vertex that is both the initial point of an edge and the end point of another edge, the value of $u_e$ is the same. Similarly, the spaces $C^k(\mathcal G)$, with $k\geq 1,$ can be defined. Furthermore, we consider the space
\[\mathfrak{D}:=\{u\in C(\mathcal G)\,:\, u_e\in C^2(I_e)\cap C^1(\bar I_e), u_e''\in L^\infty(I_e)  \text{ for any } e\in E\}.\]
Note that if $u\in \mathfrak{D}$, in principle $u\not\in C^1(\mathcal G)$; in fact, if $v\in \overline{e_1}\cap \overline{e_2}$ for some $e_1, e_2\in E$, it could happen that $u'_{e_1}(v)\neq u'_{e_2}(v).$

\subsection{Laplacian and Kirchoff condition}
Let $u\in \mathfrak D.$ For any $e\in E$, $x\in (0, l(e)),$ the Laplacian is defined as follows 
\[\Delta u(x)=u''(x)\,.\]
Furthermore, we denote the outer normal derivative of $u_e$ at the vertex $x\in \mathcal V$, by $\frac{du_e}{dn}(x)$.  Clearly,
\[\frac{d u_e(x)}{dn}=\begin{cases}
u'_e(x) & \text{ if } j(e)=x\\
-u'_e(x) & \text{ if } i(e)=x\,.
\end{cases}\]
We also define, for any $x\in \mathcal V$,
\begin{equation}\label{K}
[\mathcal K(u)](x):=\sum_{e\ni v}\frac{du_e}{dn}(x)\,.
\end{equation}
We shall impose that for any $x\in \mathcal V$
\begin{equation}\label{e502}
[\mathcal K(u)](x)=0,
\end{equation}
namely, the {\em Kirchhoff transmission condition} holds. Observe that whenever $x\in \mathcal V\cap \mathcal G$, \eqref{e502} becomes the usual
{\em Neumann boundary condition}.  

Let $x\in \mathcal V$; consider the edges $e_k\in E$, with $k=1, \ldots, \operatorname{deg}_x^+$, for which $j(e_k)=x,$ and the edges $f_k\in E$, with $k=1, \ldots,  \operatorname{deg}_x^-$, for which $i(f_k)=x.$
Note that 
\begin{equation}\label{e502u}
[\mathcal K(u)](x)=\sum_{k=1}^{ \operatorname{deg}_x^+}u'_{e_k}(x)-\sum_{k=1}^{ \operatorname{deg}_x^-}u'_{f_k}(x)\,.
\end{equation}

\subsection{Definition of solutions}
\begin{definition}\label{defsol1}
We say that $u\in \mathfrak D$ is a solution of \eqref{problema} whenever, for each $e\in E$,
\begin{equation}\label{e501}
u_e''(x) - V(x) u_e(x) = 0 \quad \text{ for all } x\in (0, l(e)),
\end{equation}
and, for each $x\in \mathcal V$, \eqref{e502} is fulfilled.
\end{definition}

Observe that for functions $f:\mathcal G\times (0, T)\to \mathbb R$, we use the notation $(\cdot)'$ when referring to the differentiation w.r.t. the variable $x$. On the other hand, the differentiation w.r.t. the time variable $t$ is denoted by $\partial_t(\cdot)$.

\begin{definition}\label{defsol2}
We say that $u\in C(\mathcal G\times (0, T])$ is a solution of \eqref{problema2} whenever 
\begin{itemize}
\item for each $x\in \mathcal G,$ $t\mapsto u(x, t)$ is of class $C^1$ in $(0, T]$, 
\item  for each $t\in (0, T]$, $u(\cdot, t)\in \mathfrak{D}$,
\item for all $e\in E$, \begin{equation}\label{e503}
\rho(x) \partial_t u_e(x,t) - u_e''(x) = 0 \quad \text{ for all } x\in (0, l(e)),
\end{equation}
\item for all $v\in\mathcal V$, 
\begin{equation}\label{e504}
[\mathcal K(u)](v, t)=0 \quad \text{ for all } \, t\in (0, T]\,.
\end{equation}
\end{itemize}
Furthermore, we say that $u\in C(\mathcal G\times [0, T])$ is a solution of the Cauchy problem \eqref{problema2}-\eqref{incond} whenever it is a solution of \eqref{problema2} and $u(x, 0)=0$ for any $x\in \mathcal G.$
\end{definition}

\section{Statement of the main results}\label{statement}\setcounter{equation}{0}

Let $x_0\in\G$ be a given reference point. For any $\beta>0$ we define the function
\begin{equation}\label{weight}
\varphi_\beta(x):=e^{-\beta d(x,x_0)}\quad \text{for any}\,\,\,x\in\G\,.
\end{equation}
Furthermore, for any $\lambda>0$ and $k\ge1$ we define the function
\begin{equation}\label{weight2}
\varphi_\lambda(x):=[d(x,x_0)+k]^{-\lambda}\quad \text{for any}\,\,\,x\in\G\,.
\end{equation}

In what follows, we denote by $$B_R(x_0):=\{x\in\G \,:\,d(x,x_0)<R\}$$ the ball of radius $R>0$ centered at $x_0$. We will also say that each vertices of the graph has more exiting edges than the entering ones if the inbound degree is smaller than the outbound degree, i.e.
\begin{equation}\label{degAss}
\operatorname{deg}_v^+\le \operatorname{deg}_v^-\quad\text{for all}\,\,\,v\in \mathcal V\,.
\end{equation}

\subsection{Elliptic equation}
In this subsection we collect the uniqueness results concerning equation \eqref{problema}, i.e. our elliptic equation.

\begin{theorem}\label{teo1}
Let hypotheses \eqref{Vbound} and \eqref{ipotesigrafo} be satisfied. Let $u$ be a solution to problem \eqref{problema} such that
\begin{equation}\label{growthcond1}
\int_{B_R(x_0)}|u(x)|^p\,d\mu(x)\le C e^{\beta R}\,,
\end{equation}
for some $C>0$, for all $R>1$. Furthermore, assume that
\begin{equation}\label{h1}
0<\beta < \frac{\sqrt{p\,V_0}}{2} \,,
\end{equation}
and suppose either one of the following conditions hold:
\begin{itemize}
\item[(i)] $p\geq 2$; 
\item[(ii)] $1\le p< 2$ and \eqref{degAss} holds.
\end{itemize}
Then
$$u(x) = 0 \quad \text{ for any } \; x\in \G.$$
\end{theorem}

\begin{theorem}\label{teo3}
Let hypotheses \eqref{Vbound2} and \eqref{ipotesigrafo} be satisfied. Let $u$ be a solution to problem \eqref{problema} such that
\begin{equation}\label{growthcond2}
\int_{B_R(x_0)}|u(x)|^p\,d\mu(x)\le \tilde C (R+k)^{\lambda}\,,
\end{equation}
for some $\tilde C>0$, for all $R>1$. Furthermore, assume either one of the following conditions
\begin{itemize}
\item[(i)] $p\geq 2$,  
\begin{equation}\label{h3}
0<\lambda<\sqrt{pV_0}\,,
\end{equation}
\item[(ii)] \eqref{degAss} holds, $1\le p< 2$,  
\begin{equation}\label{h4}
0<\lambda<\frac{\sqrt{1+4pV_0}-1}{2}\,.
\end{equation}
\end{itemize}
Then
$$u(x) = 0 \quad \text{ for any } \; x\in \G.$$
\end{theorem}

We conclude the subsection with the following

\begin{remark}\label{rem1}
\begin{itemize}
\item[(i)] If $u\in L^p_{\varphi_\beta}(\G)$ with $\varphi_\beta$ as in \eqref{weight}, then for some $C>0$ we can write, for all $R>1$,
$$
C\ge\int_{\G}|u(x)|^p\varphi_\beta(x)\,d\mu(x)\ge\int_{B_R(x_0)}|u(x)|^pe^{-\beta d(x,x_0)}\,d\mu(x)\ge e^{-\beta R}\int_{B_R(x_0)}|u(x)|^p\,d\mu(x)\,.
$$
Thus in particular one has \eqref{growthcond1}.
\item[(ii)] If $u\in L^p_{\varphi_\lambda}(\G)$ with $\varphi_\beta$ as in \eqref{weight2}, then for some $\tilde C>0$ we can write, for all $R>1$,
$$
\tilde C\ge\int_{\G}|u(x)|^p\varphi_\lambda(x)\,d\mu(x)\ge\int_{B_R(x_0)}|u(x)|^p[d(x,x_0)+k]^{-\lambda}\,d\mu(x)\ge ( R+k)^{-\lambda}\int_{B_R(x_0)}|u(x)|^p\,d\mu(x)\,.
$$
Thus in particular one has \eqref{growthcond2}.
\end{itemize}
\end{remark}

\subsection{Parabolic equation}

For $\varphi_\beta$ and $\varphi_\lambda$ we keep the notation above.

\begin{theorem}\label{teo5}
Let hypotheses \eqref{rhobound} and \eqref{ipotesigrafo} be satisfied. Let $u$ be a solution to problem \eqref{problema2}-\eqref{incond} such that 
\begin{equation}\label{growthcond3}
\int_0^\tau\int_{B_R(x_0)}|u(x,t)|^p\,d\mu(x)dt\le C e^{\beta R}\,,
\end{equation}
for some $C>0$, for all $R>1$ and for all $\tau\in(0,T]$. Furthermore, assume either one of the following
\begin{itemize}
\item[(i)] $p\geq 2$ and $\beta>0$;
\item[(ii)] \eqref{degAss} holds, $1\le p<2$, and $\beta>0$.
\end{itemize}
Then
$$u(x,t) = 0 \quad \text{ for any } \; x\in \G,\,\,t\in(0,T].$$
\end{theorem}


\begin{theorem}\label{teo7}
Let hypotheses  \eqref{rhobound2} and \eqref{ipotesigrafo} be satisfied. Let $u$ be a solution to problem \eqref{problema2}-\eqref{incond} such that  
\begin{equation}\label{growthcond4}
\int_0^\tau\int_{B_R(x_0)}|u(x,t)|^p\,d\mu(x)dt\le \tilde C ( R+k)^{\lambda}\,,
\end{equation}
for some $\tilde C>0,$ for all $R>1$ and for all $\tau\in(0,T]$. Assume either one of the following
\begin{itemize}
\item[(i)]  $p\geq 2$ and $\lambda>0$;
\item[(ii)] \eqref{degAss} holds, $1\le p<2$, and $\lambda>0$.
\end{itemize}
Then
$$u(x,t) = 0 \quad \text{ for any } \; x\in \G,\,\,t\in(0,T].$$
\end{theorem}


Similarly to Remark \ref{rem1}, we can observe the following

\begin{remark}\label{rem2}
\begin{itemize}
\item[(i)] If $u\in L^p_{\varphi_\beta}(\G\times(0,T))$ with $\varphi_\beta$ as in \eqref{weight}, then for every $\tau\in(0,T)$ and for some $C>0$ we can write, for all $R>1$,
$$
\begin{aligned}
C\ge\int_0^\tau\int_{\G}|u(x,t)|^p&\varphi_\beta(x)\,d\mu(x)dt\ge\int_0^\tau\int_{B_R(x_0)}|u(x,t)|^pe^{-\beta d(x,x_0)}\,d\mu(x)dt\\
&\ge e^{-\beta R}\int_0^\tau\int_{B_R(x_0)}|u(x,t)|^p\,d\mu(x)dt\,.
\end{aligned}
$$
Thus in particular one has \eqref{growthcond3}.
\item[(ii)] If $u\in L^p_{\varphi_\lambda}(\G\times(0,T))$ with $\varphi_\lambda$ as in \eqref{weight2}, then for every $\tau\in(0,T)$ and for some $\tilde C>0$ we can write, for all $R>1$,
$$
\begin{aligned}
\tilde C\ge\int_0^\tau\int_{\G}|u(x,t)|^p&\varphi_\lambda(x)\,d\mu(x)dt\ge\int_0^\tau\int_{B_R(x_0)}|u(x,t)|^p[d(x,x_0)+k]^{-\lambda}\,d\mu(x)dt\\
&\ge [ R+k]^{-\lambda}\int_0^\tau\int_{B_R(x_0)}|u(x,t)|^p\,d\mu(x)dt\,.
\end{aligned}
$$
Thus in particular one has \eqref{growthcond4}.
\end{itemize}
\end{remark}

\section{Auxiliary results} \label{sec2}\setcounter{equation}{0}

\subsection{Elliptic equation}
Before delving into the details of the proofs of the main results, we need to establish some preliminary facts that will play a key role in the next section.

\begin{proposition}\label{prop41}
Let $u$ be a solution to equation \eqref{problema}. Let $\eta, \xi:\mathcal G\to\R$ be such that
\begin{itemize}
\item $\eta\ge0, \operatorname{supp} \eta$ is compact, $\eta\in C^1(\mathcal G);$
\item $\xi \in C^1(\mathcal G)$.
\end{itemize}
Then, for any $p\ge 2$,
\begin{equation*}\label{eq42}
\begin{aligned}
\int_{\G}&|u(x)|^p\eta^2(x)e^{\xi(x)}\left[p\,V(x)-\left(\xi'(x)\right)^2\right]\,d\mu(x)\le4\int_{\G}|u(x)|^p\left(\eta'(x)\right)^2e^{\xi(x)}\,d\mu(x)\\
\end{aligned}
\end{equation*}
\end{proposition}

\begin{proof}
For any $\alpha>0$, we define
\begin{equation}\label{eq43}
\pi_{\alpha}(u):=(u^2+\alpha)^{\frac p4}\,.
\end{equation}
Now, for any $x\in \G\setminus\mathcal V$, we consider the expression
\begin{equation}\label{eq44}
\Delta \pi_{\alpha}(u(x))-V(x)\,\pi_{\alpha}(u(x)).
\end{equation}
We multiply \eqref{eq44} by $\pi_{\alpha}(u(x))\,\eta^2(x)e^{\xi(x)}$ and then we integrate over $\G$. Thus we take into account
\begin{equation}\label{eq45}
\int_{\G}\Delta \pi_{\alpha}(u(x))\pi_{\alpha}(u(x))\,\eta^2(x)e^{\xi(x)}d\mu(x)\,-\int_{\G}\,V(x)\,\pi^2_{\alpha}(u(x))\,\eta^2(x)e^{\xi(x)}d\mu(x)\,.
\end{equation}
Set
\begin{equation}\label{eq46}
I:=\int_{\G}\Delta \pi_{\alpha}(u(x))\pi_{\alpha}(u(x))\,\eta^2(x){e^\xi(x)}\,d\mu(x)\,.
\end{equation}
Then, by \eqref{int1} and integrating by parts we get
\begin{equation}\label{eq47}
\begin{aligned}
I&= \sum_{e\in E}\int_0^{l_e}[\pi_{\alpha}(u_e(x))]{''}\pi_{\alpha}(u_e(x))\,\eta_e^2(x){e^{\xi_e(x)}}\,dx\\
&=-\sum_{e\in E}\int_0^{l_e} \pi_{\alpha}'(u_e(x)) \left[\pi_{\alpha}(u_e(x))\eta_e^2(x)e^{\xi_e(x)}\right]'\, dx\\
&\quad+\sum_{x\in V}[\mathcal{K}(\pi_{\alpha}(u))](x) \pi_{\alpha}(u(x))\eta^2(x)e^{\xi(x)}\\
&=-\sum_{e\in E}\int_0^{l_e} (\pi_{\alpha}'(u_e(x)))^2 \eta_e^2(x)e^{\xi_e(x)}\, dx \\
&\quad -2\sum_{e\in E}\int_0^{l_e} \pi_{\alpha}'(u_e(x))\pi_{\alpha}(u_e(x)) \eta_e(x)\eta_e'(x)e^{\xi_e(x)}\, dx\\
&\quad-\sum_{e\in E}\int_0^{l_e} \pi_{\alpha}'(u_e(x))\pi_{\alpha}(u_e(x)) \eta_e^2(x)\xi_e'(x)e^{\xi_e(x)}\, dx\\
&\quad+\sum_{x\in V}[\mathcal{K}(u)](x)\pi_{\alpha}'(u(x)) \pi_{\alpha}(u(x))\eta^2(x)e^{\xi(x)}\\
&=:\sum_{e\in E}\{J_1+J_2+J_3\}\,,
\end{aligned}
\end{equation}
where in the last equality we used \eqref{e502}. By Young's inequality with exponent $2$, we have, for any $\delta_1>0$ and $\delta_2>0$,
\begin{equation}\label{eq49}
\begin{aligned}
J_2&=-2\int_0^{l_e} \pi_{\alpha}'(u_e(x))\pi_{\alpha}(u_e(x)) \eta_e(x)\eta'_e(x)e^{\xi_e(x)}\, dx\\
&\le\delta_1\int_0^{l_e} (\pi_{\alpha}'(u_e(x)))^2 \eta_e^2(x)e^{\xi_e(x)}\, dx+\frac{1}{\delta_1}\int_0^{l_e} \pi_{\alpha}^2(u_e(x))(\eta_e'(x))^2e^{\xi_e(x)}\, dx\,;
\end{aligned}
\end{equation}
and
\begin{equation}\label{eq412}
\begin{aligned}
J_3&=-\int_0^{l_e} \pi_{\alpha}'(u_e(x))\pi_{\alpha}(u_e(x)) \eta_e^2(x)\xi'_e(x)e^{\xi_e(x)}\, dx\\
&\le\frac{\delta_2}{2}\int_0^{l_e} (\pi_{\alpha}'(u_e(x)))^2 \eta_e^2(x)e^{\xi_e(x)}\,dx+\frac{1}{2\delta_2}\int_0^{l_e} \pi_{\alpha}^2(u_e(x))\eta_e^2(x)(\xi_e'(x))^2e^{\xi_e(x)}\, dx\,.
\end{aligned}
\end{equation}
By combining together \eqref{eq46}, \eqref{eq47}, \eqref{eq49} and \eqref{eq412} we get
\begin{equation*}
\begin{aligned}
\int_{\G}&\Delta \pi_{\alpha}(u(x))\pi_{\alpha}(u(x))\,\eta^2(x){e^\xi(x)}\,d\mu(x)\\
&\le \sum_{e\in E}\left\{\left(-1+\delta_1+\frac{\delta_2}{2}\right)\int_0^{l_e} (\pi_{\alpha}'(u_e(x)))^2 \eta_e^2(x)e^{\xi_e(x)}\, dx\right.\\
&\left. +\frac{1}{\delta_1}\int_0^{l_e} \pi_{\alpha}^2(u_e(x))(\eta_e'(x))^2e^{\xi_e(x)}\, dx+\frac{1}{2\delta_2}\int_0^{l_e} \pi_{\alpha}^2(u_e(x))\eta_e^2(x)(\xi_e'(x))^2e^{\xi_e(x)}\, dx\right\}\,.
\end{aligned}
\end{equation*}
We now choose $\delta_1=\frac 12$ and $\delta_2=1$ in such a way that $-1+ \delta_1+ \frac{\delta_2}{2}=0$. Consequently, the latter yields
\begin{equation}\label{eq414}
\begin{aligned}
\int_{\G}&\Delta \pi_{\alpha}(u(x))\pi_{\alpha}(u(x))\,\eta^2(x){e^\xi(x)}\,d\mu(x)\\
&\le\sum_{e\in E}\left\{2\int_0^{l_e} \pi_{\alpha}^2(u_e(x))(\eta_e'(x))^2e^{\xi_e(x)}\, dx+\frac12\int_0^{l_e} \pi_{\alpha}^2(u_e(x))\eta_e^2(x)(\xi'_e(x))^2e^{\xi_e(x)}\, dx\right\}\,.
\end{aligned}
\end{equation}
By substituting \eqref{eq414} into \eqref{eq45}, we get
\begin{equation}\label{eq416}
\begin{aligned}
\int_{\G}\Delta \pi_{\alpha}(u(x))&\pi_{\alpha}(u(x))\,\eta^2(x)e^{\xi(x)}d\mu(x)\,-\int_{\G}\,V(x)\,\pi^2_{\alpha}(u(x))\,\eta^2(x)e^{\xi(x)}d\mu(x)\\
&\le2\sum_{e\in E}\int_0^{l_e} \pi_{\alpha}^2(u_e(x))(\eta_e'(x))^2e^{\xi_e(x)}\, dx\\
&\quad+\sum_{e\in E}\int_0^{l_e} \pi_{\alpha}^2(u_e(x))\eta_e^2(x)\left[\frac12(\xi_e'(x))^2-V_e(x)\right]e^{\xi_e(x)}\, dx\,.
\end{aligned}
\end{equation}

Moreover, observe that, due to the assumption $p\ge2$, we get, for all $x\in \G\setminus\mathcal V$
\begin{equation}\label{eq417}
\begin{aligned}
\left(\Delta \pi_{\alpha}(u(x))\right)_e&=\pi''_\alpha(u_e(x))=\left[\frac p2 u_e(x)(u_e^2+\alpha)^{p/4-1}u_e'(x)\right]\\
&=\frac p2(u_e^2+\alpha)^{p/4-1}(u_e'(x))^2+\frac p2(u_e^2+\alpha)^{p/4-1}u_e(x)u_e''(x)\\
&\quad+p\left(\frac p4-1\right)u_e^2(x)(u_e'(x))^2(u_e^2+\alpha)^{p/4-2}\\
&=\frac p4(u_e^2+\alpha)^{p/4-2}(u_e')^2\left\{\frac{p-2}{2}u_e^2(x)+\alpha\right\}+\frac p2(u_e^2+\alpha)^{p/4-1}u_e(x)u_e''(x)\\
&\ge\frac p2\left(u_e^2+\alpha\right)^{p/4-1}u_e(x)u_e''(x)=\pi_{\alpha}'(u_e(x))u_e''(x)\,.
\end{aligned}
\end{equation}
Thus, the left-hand-side of \eqref{eq416}, by \eqref{eq417} and \eqref{problema} satisfies the inequality
\begin{equation}\label{eq418}
\begin{aligned}
\int_{\G}&\left[\Delta \pi_{\alpha}(u(x))-V(x)\,\pi_{\alpha}(u(x))\right]\pi_{\alpha}(u(x))\,\eta^2(x)e^{\xi(x)}d\mu(x)\,\\
&\ge \sum_{e\in E}\int_0^{l_e}\left[\pi'_\alpha(u_e(x)) u_e''(x)-V_e(x)\pi_{\alpha}(u_e(x)) \right]\pi_{\alpha}(u_e(x))\eta_e^2(x){e^{\xi_e(x)}}\,dx\\
&=\sum_{e\in E}\int_0^{l_e}V_e(x)\left[\pi'_\alpha(u_e(x))u_e(x)-\pi_{\alpha}(u_e(x)) \right]\pi_{\alpha}(u_e(x))\eta_e^2(x){e^{\xi_e(x)}}\,dx.
\end{aligned}
\end{equation}
By combining together \eqref{eq416} and \eqref{eq418} we get
\begin{equation}\label{eq421}
\begin{aligned}
\sum_{e\in E}\int_0^{l_e}&V_e(x)\pi'_\alpha(u_e(x))u_e(x)\pi_{\alpha}(u_e(x))\eta_e^2(x){e^{\xi_e(x)}}\,dx\\
&\le2\sum_{e\in E}\int_0^{l_e} \pi_{\alpha}^2(u_e(x))(\eta_e'(x))^2e^{\xi_e(x)}\, dx\\
&\quad+\sum_{e\in E}\frac12\int_0^{l_e} (\xi_e'(x))^2\pi_{\alpha}^2(u_e(x))\eta_e^2(x)e^{\xi_e(x)}\, dx\,.
\end{aligned}
\end{equation}
Now, since
\[\pi_\alpha(u)\to |u|^{\frac p2}, \quad \text{ as } \alpha\to 0^+,\]
letting $\alpha\to 0^+$ in \eqref{eq421}, we arrive to
\begin{equation*}\label{eq312}
\begin{aligned}
\frac12\sum_{e\in E}\int_0^{l_e}&|u_e(x)|^p\eta_e^2(x){e^{\xi_e(x)}}\left[pV_e(x)-(\xi_e'(x))^2\right]\,dx\le2\sum_{e\in E}\int_0^{l_e} |u_e(x)|^p(\eta_e'(x))^2e^{\xi_e(x)}\, dx\,,
\end{aligned}
\end{equation*}
which yields the thesis.
\end{proof}

\begin{corollary}\label{cor41}
Let $u$ be a solution to equation \eqref{problema}. Let $\eta, \zeta:\mathcal G\to\R$ be such that
\begin{itemize}
\item $\eta\ge0$, $\operatorname{supp} \eta$ is compact, $\eta\in C^1(\mathcal G)$;
\item $\zeta\in C^1(\mathcal G), \zeta>0$.
\end{itemize}
Then, for any $p\ge 2$,
\begin{equation*}\label{eq42}
\begin{aligned}
\int_{\G}&|u(x)|^p\eta^2(x){\zeta(x)}\left[p\,V(x)-\left(\frac{\zeta'(x)}{\zeta(x)}\right)^2\right]\,d\mu(x)\le4\int_{\G}|u(x)|^p\left(\eta'(x)\right)^2{\zeta(x)}\,d\mu(x)\\
\end{aligned}
\end{equation*}
\end{corollary}

\begin{proof}
The thesis follows by applying Proposition \ref{prop41} with $\xi_e(x)=\log\zeta_e(x)$ for every $e\in E$.
\end{proof}

\begin{proposition}\label{prop42}
Let $u$ be a solution to problem \eqref{problema}. Let $p\ge1$ and $f\in C^2(\mathcal G)$ be a nonnegative function with compact support. Then
\begin{equation}\label{eqlemma}
\int_{\G}|u(x)|^p\left[-\Delta f(x)+pV(x) f(x)\right] \,d\mu(x)\le-\sum_{x\in V}|u(x)|^p[\mathcal K(f)](x)\,.
\end{equation}
\end{proposition}

\begin{proof}
For any $\alpha>0$, we define
\begin{equation}\label{G}
\phi_{\alpha}(u):=(u^2+\alpha)^{\frac p2}\,.
\end{equation}
We consider
\begin{equation}\label{start}
\int_{\G}\left[\Delta \phi_{\alpha}(u(x))-V(x)\,\phi_{\alpha}(u(x))\right]f(x)\,d\mu(x).
\end{equation}
By \eqref{int1}, integrating by parts twice, we obtain
\begin{equation}\label{eq33b}
\begin{aligned}
\int_{\G}&\left[\Delta \phi_{\alpha}(u(x))-V(x)\,\phi_{\alpha}(u(x))\right]f(x)\,d\mu(x)\\
&=\sum_{e\in E}\left\{\int_0^{l_e}\left[\phi_\alpha''(u_e(x))-V_e(x)\phi_\alpha(u_e(x))\right]f_e(x)\,dx\right\}\\
&=-\sum_{e\in E}\int_0^{l_e}\phi_\alpha'(u_e(x))f_e'(x)\,dx+\sum_{x\in V}f(x)\left[\mathcal K(\phi_\alpha(u))\right](x)\\
&\quad\quad-\sum_{e\in E}\int_0^{l_e}V_e(x)\phi_\alpha(u_e(x))f_e(x)\,dx\\
&=\sum_{e\in E}\int_0^{l_e}\phi_\alpha(u_e(x))f_e''(x)\,dx-\sum_{x\in V}\phi_\alpha(u(x))\left[\mathcal K(f)\right](x)\\
&\quad\quad+\sum_{x\in V}f(x)\left[\mathcal K(\phi_\alpha(u))\right](x)-\sum_{e\in E}\int_0^{l_e}V_e(x)\phi_\alpha(u_e(x))f_e(x)\,dx\\
\end{aligned}
\end{equation}
By \eqref{e502},
\begin{equation}\label{e33g}
\sum_{x\in V}f(x)\left[\mathcal K(\phi_\alpha(u))\right](x)=\sum_{x\in V}f(x)\phi_\alpha'(u)(x)\left[\mathcal K(u)\right](x)=0.
\end{equation}
Moreover, similarly to \eqref{eq417}, by using that $p\geq 1$, we get, for all $x\in\G\setminus\mathcal V$
\begin{equation}\label{e33f}
\Delta\phi_\alpha(u_e(x))=\phi_\alpha''(u_e(x))\geq \phi'_\alpha(u_e(x)) u_e''(x)\,.
\end{equation}
Therefore, \eqref{eq33b}, \eqref{e33g} and \eqref{e33f} together with problem \eqref{problema} yield
\begin{equation}\label{e34f}
\begin{aligned}
\sum_{e\in E}&\left\{\int_0^{l_e}\left[\phi_\alpha'(u_e(x))u_e(x)-\phi_\alpha(u_e(x))\right]V_e(x)f_e(x)\,dx\right\}\\
&=\sum_{e\in E}\left\{\int_0^{l_e}\left[\phi_\alpha'(u_e(x))u_e''(x)-V_e(x)\phi_\alpha(u_e(x))\right]f_e(x)\,dx\right\}\\
&\le \sum_{e\in E}\left\{\int_0^{l_e}\phi_\alpha(u_e(x))\left[f_e''(x)-V_e(x)f_e(x)\right]\,dx\right\}-\sum_{x\in V}\phi_\alpha(u(x))\left[\mathcal K(f)\right](x)\\
\end{aligned}
\end{equation}
Letting $\alpha\to 0^+$ in \eqref{e34f},
 $$
 \phi_\alpha(u_e(x))\to |u_e(x)|^p,\quad \phi_\alpha'(u_e(x))\to p\, |u_e(x)|^{p-1}
 $$
 hence we get
 $$
\sum_{e\in E}\int_0^{l_e}|u_e(x)|^p\left[-f_e''(x)+pV_e(x)f_e(x)\right]dx\le-\sum_{x\in V}|u(x)|^p\left[\mathcal K(f)\right](x)\,,
 $$
 which yield the thesis.
 \end{proof}

\subsection{Parabolic equation}

\begin{proposition}\label{prop51}
Let $u$ be a solution to equation \eqref{problema2} with initial condition \eqref{incond}. Let $\eta:\mathcal G\to\R$ and $\omega:\G\times[0,T]\to\R$ be such that
\begin{itemize}
\item $\eta\ge0, \operatorname{supp} \eta$ is compact, $\eta\in C^1(\mathcal G)$;
\item $\omega\in C^{2,1}_{x,t}(\mathcal G\times [0, T]$.
\end{itemize}
Then, for any $p\ge 2$ and any $\tau\in(0,T]$,
\[
\begin{aligned}
\int_{\G}\,\rho(x)\,&|u(x,\tau)|^pe^{\omega(x,\tau)}\eta^2(x)\,d\mu(x)\le 4\int_0^\tau\int_{\G}|u(x,t)|^p(\eta'(x))^2e^{\omega(x,t)}\, d\mu(x)\,dt\\
&\quad+\int_0^\tau \int_{\G}|u(x,t)|^p\left[\rho(x)\,\partial_t \omega(x,t)+(\omega'(x,t))^2\right]e^{\omega(x,t)}\eta^2(x)\,d\mu(x)\,dt\,.
\end{aligned}
\]
\end{proposition}

\begin{proof}
For any $\alpha>0$, we define
\begin{equation}\label{eq53}
\pi_{\alpha}(u):=(u^2+\alpha)^{\frac p4}\,.
\end{equation}
Now, for any $x\in \G\setminus \mathcal V$ and $t\in(0,\tau]$, we consider the expression
\begin{equation}\label{eq54}
\rho(x)\,\partial_t\pi_{\alpha}(u(x,t))-\Delta \pi_{\alpha}(u(x,t)).
\end{equation}
We multiply \eqref{eq54} by $\pi_{\alpha}(u(x,t))\,\eta^2(x)e^{\omega(x,t)}$ and then we integrate over $\G$. Thus
\begin{equation}\label{eq55}
\begin{aligned}
\int_{\G} \rho(x)\,\partial_t \pi_{\alpha}(u(x,t))&\pi_{\alpha}(u(x,t))\,\eta^2(x)e^{\omega(x,t)}d\mu(x)\,\\
&-\int_{\G}\,\Delta \pi_{\alpha}(u(x,t))\pi_{\alpha}(u(x,t))\,\eta^2(x)e^{\omega(x,t)}d\mu(x)\,.
\end{aligned}
\end{equation}
Set
\begin{equation}\label{eq56}
\begin{aligned}
I&:=\int_{\G}\Delta \pi_{\alpha}(u(x,t))\pi_{\alpha}(u(x,t))\,\eta^2(x){e^\omega(x,t)}\,d\mu(x)\,\\
J&:=\int_{\G} \rho(x)\,\partial_t \pi_{\alpha}(u(x,t))\pi_{\alpha}(u(x,t))\,\eta^2(x)e^{\omega(x,t)}d\mu(x)
\end{aligned}
\end{equation}
We perform on $I$ the same computation done in the proof of Proposition \ref{prop41} on \eqref{eq46}; this yields to
$$
\begin{aligned}
I&\le\frac12\int_{\G} \pi_{\alpha}^2(u(x,t))\eta^2(x)(\omega'(x,t))^2e^{\omega(x,t)}\, d\mu(x)+2\int_{\G} \pi_{\alpha}^2(u(x,t))(\eta'(x))^2e^{\omega(x,t)}\, d\mu(x)
\end{aligned}
$$
Moreover, we can write $J$ as follows
$$
\begin{aligned}
J&=\frac12\int_{\G} \rho(x)\partial_t \left[\pi_{\alpha}^2(u(x,t))\right]\,\eta^2(x)e^{\omega(x,t)}d\mu(x)\\
&=\frac12\int_{\G} \rho(x)\partial_t \left[\pi_{\alpha}^2(u(x,t))e^{\omega(x,t)}\right]\eta^2(x)d\mu(x)-\frac12\int_{\G} \rho(x)\pi_{\alpha}^2(u(x,t))e^{\omega(x,t)}\partial_t(\omega(x,t))\eta^2(x)d\mu(x)\,.
\end{aligned}
$$
Combining the estimates on $I$ and $J$ together with \eqref{eq55} and \eqref{eq56}, we get
\begin{equation}\label{eq516}
\begin{aligned}
\int_{\G} \rho(x)&\partial_t \pi_{\alpha}(u(x,t))\pi_{\alpha}(u(x,t))\,\eta^2(x)e^{\omega(x,t)}d\mu(x)\\
&-\int_{\G}\,\Delta \pi_{\alpha}(u(x,t))\pi_{\alpha}(u(x,t))\,\eta^2(x)e^{\omega(x,t)}d\mu(x)\\
&\ge\frac12\int_{\G} \rho(x)\,\partial_t \left[\pi_{\alpha}^2(u(x,t))e^{\omega(x,t)}\right]\,\eta^2(x)d\mu(x)\\
&\quad-\frac12\int_{\G} \rho(x)\,\pi_{\alpha}^2(u(x,t))e^{\omega(x,t)}\left[\partial_t \left({\xi(x,t)}\right)+(\omega'(x))^2\right]\,\eta^2(x)d\mu(x)\\
&\quad-2\int_{\G} \pi_{\alpha}^2(u(x,t))(\eta'(x))^2e^{\omega(x,t)}\, d\mu(x)\,.\\
\end{aligned}
\end{equation}
Furthermore, we may estimate $I$ defined in \eqref{eq56} from below, in fact, since $p\ge2$, we may use inequality \eqref{eq417} and therefore write
\begin{equation}\label{eq517}
I\ge\int_{\G}\pi_{\alpha}'(u(x,t))\Delta u(x,t)\pi_{\alpha}(u(x,t))\,\eta^2(x){e^ \omega(x,t)}\,d\mu(x)\,.
\end{equation}
Thus, the left-hand-side of \eqref{eq516}, by \eqref{eq517} and \eqref{problema2} satisfies the inequality
\begin{equation}\label{eq518}
\begin{aligned}
\int_{\G} \rho(x)\partial_t& \pi_{\alpha}(u(x,t))\pi_{\alpha}(u(x,t))\eta^2(x)e^{\omega(x,t)}d\mu(x)\\
&\quad-\int_{\G}\Delta \pi_{\alpha}(u(x,t))\pi_{\alpha}(u(x,t))\,\eta^2(x)e^{\omega(x,t)}d\mu(x)\\
&\le \int_{\G}\pi'_\alpha(u(x,t))\left[ \rho(x)\partial_tu(x,t)-\Delta u(x,t) \right]\pi_{\alpha}(u(x,t))\eta^2(x){e^{\omega(x,t)}}\,d\mu(x)=0
\end{aligned}
\end{equation}
By combining together \eqref{eq516} and \eqref{eq518} we get
\begin{equation}\label{eq521}
\begin{aligned}
\frac{\partial}{\partial t} &\left(\int_{\G} \rho(x)\pi_{\alpha}^2(u(x,t))e^{\omega(x,t)}\,\eta^2(x)d\mu(x)\right)\\
&\le\int_{\G} \pi_{\alpha}^2(u(x,t))e^{\omega(x,t)}\left[\rho(x)\partial_t \left({\omega(x,t)}\right)+(\omega'(x))^2\right]\,\eta^2(x)d\mu(x)\\
&\quad+4\int_{\G} \pi_{\alpha}^2(u(x,t))(\eta'(x))^2e^{\omega(x,t)}\, d\mu(x)\,.\\
\end{aligned}
\end{equation}
We now integrate the latter with respect to the variable $t$ in the interval $[0,\tau]$, for any $\tau\in(0,T]$, thus we write
\begin{equation}\label{eq522}
\begin{aligned}
\int_{\G}\,\rho(x)\,&\pi_{\alpha}^2(u(x,\tau))e^{\omega(x,\tau)}\eta^2(x)\,d\mu(x)-\int_{\G}\,\rho(x)\alpha^{p/2}e^{\omega(x,0)}\eta^2(x)\,d\mu(x)\\
&\le\int_0^\tau \int_{\G}\,\left[\rho(x)\,\partial_t \omega(x,t)+(\omega'(x,t))^2\right]\pi^2_{\alpha}(u(x,t))e^{\omega(x,t)}\eta^2(x)\,d\mu(x)\,dt\\
&\quad + 4\int_0^\tau\int_{\G}\pi_{\alpha}^2(u(x,t))(\eta'(x))^2e^{\omega(x,t)}\, d\mu(x)\,dt\,.
\end{aligned}
\end{equation}
where we have used that $u(x,0)=0$ by \eqref{problema2}.
Now, since
\[\pi_\alpha(u)\to |u|^{\frac p2}, \quad \text{ as } \alpha\to 0^+,\]
letting $\alpha\to 0^+$ in \eqref{eq522}, we arrive to
\begin{equation*}
\begin{aligned}
\int_{\G}\,\rho(x)\,&|u(x,\tau)|^pe^{\omega(x,\tau)}\eta^2(x)\,d\mu(x)\\
&\le\int_0^\tau \int_{\G}|u(x,t)|^p\left[\rho(x)\,\partial_t \omega(x,t)+(\omega'(x,t))^2\right]e^{\omega(x,t)}\eta^2(x)\,d\mu(x)\,dt\\
&\quad + 4\int_0^\tau\int_{\G}|u(x,t)|^p\eta'(x))^2e^{\omega(x,t)}\, d\mu(x)\,dt\,.
\end{aligned}
\end{equation*}
which yields the thesis.
\end{proof}

\begin{corollary}\label{cor51}
Let $u$ be a solution to problem \eqref{problema2}-\eqref{incond}. Let $\eta:\mathcal G\to\R$ and $\kappa:\G\times[0,T]\to\R$ be such that
\begin{itemize}
\item $\eta\ge0, \operatorname{supp} \eta$ is compact, $\eta\in C^1(\mathcal G);$
\item $\kappa>0$ and $\kappa\in C^1(\mathcal G\times[0,T])$.
\end{itemize}
Then, for any $p\ge 2$ and any $\tau\in(0,T]$,
\[
\begin{aligned}
\int_{\G}\,\rho(x)\,&|u(x,\tau)|^p{\kappa(x,\tau)}\eta^2(x)\,d\mu(x)\le 4\int_0^\tau\int_{\G}|u(x,t)|^p(\eta'(x))^2{\kappa(x,t)}\, d\mu(x)\,dt\\
&\quad+\int_0^\tau \int_{\G}|u(x,t)|^p\left[\rho(x)\,\partial_t \kappa(x,t)+\frac{(\kappa'(x,t))^2}{\kappa(x,t)}\right]\eta^2(x)\,d\mu(x)\,dt\,.
\end{aligned}
\]
\end{corollary}
\begin{proof}
The thesis follows by applying Proposition \ref{prop51} with $\omega(x,t)=\log\kappa(x,t)$ for every $x\in\G$ and $t\in(0,T]$.
\end{proof}

\begin{proposition}\label{prop52}
Let $u$ be a solution to problem \eqref{problema2}-\eqref{incond}. Let $p\ge1$ and $v\in C^{2,1}_{x,t}(\mathcal G\times [0, T])$ be a nonnegative function such that for each $t\in [0, T]$, $v(\cdot, t)$ has compact support. Then, for every $\tau\in(0,T]$,
\begin{equation}\label{eq530}
\begin{aligned}
\int_0^\tau\int_{\G}|u(x,t)|^p&\left[\rho(x)\partial_tv(x,t)+\Delta v(x,t)\right] \,d\mu(x)\,dt\\
&\ge\int_{\G}\rho(x)|u(x,\tau)|^pv(x,\tau)\,d\mu(x)+\int_0^\tau\sum_{x\in V}|u(x,t)|^p[\mathcal K(f)](x,t)\,dt\,.\\
\end{aligned}
\end{equation}
\end{proposition}

\begin{proof}
For any $\alpha>0$, we define
\begin{equation}\label{eq532}
\phi_{\alpha}(u):=(u^2+\alpha)^{\frac p2}\,.
\end{equation}
We consider
\begin{equation}\label{eq531}
\int_{\G}\left[\rho(x)\partial_t\phi_{\alpha}(u(x,t))-\Delta \phi_{\alpha}(u(x,t))\right]v(x,t)\,d\mu(x).
\end{equation}
By \eqref{int1}, integrating by parts twice, we obtain
\begin{equation}\label{eq533}
\begin{aligned}
\int_{\G}&\left[\rho(x)\partial_t\phi_{\alpha}(u(x,t))-\Delta \phi_{\alpha}(u(x,t))\right]v(x,t)d\mu(x)\\
&=\sum_{e\in E}\left\{\int_0^{l_e}\left[\rho_e(x)\partial_t\phi_{\alpha}(u_e(x,t))-\phi_{\alpha}''(u_e(x,t))\right]v_e(x,t)\,dx\right\}\\
&=-\sum_{e\in E}\int_0^{l_e}\rho_e(x)\partial_tv_e(x,t)\,\phi_{\alpha}(u_e(x,t))\,dx+\partial_t\left[\sum_{e\in E}\int_0^{l_e}\rho_e(x)v_e(x,t)\,\phi_{\alpha}(u_e(x,t))\,dx\right]\\
&\quad+\sum_{e\in E}\int_0^{l_e}\phi_\alpha'(u_e(x,t))v_e'(x,t)\,dx-\sum_{x\in V}v(x,t)\left[\mathcal K(\phi_\alpha(u))\right](x,t)\\
&=-\sum_{e\in E}\int_0^{l_e}\rho_e(x)\partial_tv_e(x,t)\,\phi_{\alpha}(u_e(x,t))\,dx+\partial_t\left[\sum_{e\in E}\int_0^{l_e}\rho_e(x)v_e(x,t)\,\phi_{\alpha}(u_e(x,t))\,dx\right]\\
&\quad-\sum_{e\in E}\int_0^{l_e}\phi_\alpha(u_e(x,t))v_e''(x,t)\,dx+\sum_{x\in V}\phi_\alpha(u(x,t))\left[\mathcal K(f)\right](x,t)\\
&\quad-\sum_{x\in V}v(x,t)\left[\mathcal K(\phi_\alpha(u))\right](x,t)\,.
\end{aligned}
\end{equation}
By \eqref{e502},
\begin{equation}\label{eq534}
\sum_{x\in V}v(x,t)\left[\mathcal K(\phi_\alpha(u))\right](x,t)=\sum_{x\in V}v(x,t)\phi_\alpha'(u)(x,t)\left[\mathcal K(u)\right](x,t)=0.
\end{equation}
Moreover, similarly to \eqref{eq417}, by using that $p\geq 1$, we get, for all $x\in\G\setminus\mathcal V$ and $t\in(0,\tau]$
\begin{equation*}
\Delta\phi_\alpha(u_e(x,t))=\phi_\alpha''(u_e(x,t))\geq \phi'_\alpha(u_e(x,t)) u_e''(x,t)\,,
\end{equation*}
and
$$
\rho_e(x)\partial_t\phi_{\alpha}(u_e(x,t))=\rho_e(x)\phi'_\alpha(u_e(x,t)) \partial_t(u_e(x,t))\,.
$$
Therefore, by using problem \eqref{problema2}, \eqref{eq533} reads
\begin{equation}\label{eq535}
\begin{aligned}
\sum_{e\in E}&\left\{\int_0^{l_e}\left[\rho_e(x)\partial_tv_e(x,t)+v_e''(x,t)\right]\phi_\alpha(u_e(x,t))\,dx\right\}\\
&\ge \partial_t\left[\sum_{e\in E}\int_0^{l_e}\rho_e(x)v_e(x,t)\,\phi_{\alpha}(u_e(x,t))\,dx\right]\\
&\quad +\sum_{x\in V}\phi_\alpha(u(x,t))\left[\mathcal K(v)\right](x,t)\,.
\end{aligned}
\end{equation}
We now integrate the \eqref{eq535} with respect to the variable $t$ in the interval $[0,\tau]$, for any $\tau\in(0,T]$, since $u(x,0)=0$ we get
$$
\begin{aligned}
\int_0^\tau\sum_{e\in E}&\left\{\int_0^{l_e}\left[\rho_e(x)\partial_tv_e(x,t)+v_e''(x,t)\right]\phi_\alpha(u_e(x,t))\,dx\right\}\\
&\ge \sum_{e\in E}\int_0^{l_e}\rho_e(x)v_e(x,\tau)\,\phi_{\alpha}(u_e(x,\tau))\,dx-\sum_{e\in E}\int_0^{l_e}\rho_e(x)v_e(x,0)\,\alpha^{p/2},dx\\
&\quad +\int_0^\tau\sum_{x\in V}\phi_\alpha(u(x,t))\left[\mathcal K(v)\right](x,t)\,.
\end{aligned}
$$
Letting $\alpha\to 0^+$,
 $$
 \phi_\alpha(u_e(x,t))\to |u_e(x,t)|^p,\quad \phi_\alpha'(u_e(x,t))\to p\, |u_e(x,t)|^{p-1}\,.
 $$
 Hence we get
 $$
\begin{aligned}
\int_0^\tau\sum_{e\in E}&\left\{\int_0^{l_e}\left[\rho_e(x)\partial_tv_e(x,t)+v_e''(x,t)\right]|u_e(x,t)|^p\,dx\right\}\\
&\ge \sum_{e\in E}\int_0^{l_e}\rho_e(x)v_e(x,\tau)\,|u_e(x,\tau)|^p\,dx +\int_0^\tau\sum_{x\in V}|u_e(x,t)|^p\left[\mathcal K(v)\right](x,t),
\end{aligned}
 $$
 which yields the thesis.
 \end{proof}

\section{Proofs the main results for the elliptic equation}\label{sec5}\setcounter{equation}{0}

\subsection{Potential bounded away from zero}

 For any $\alpha>0$, we define the function
 \begin{equation}\label{xi}
\xi(x):=-\alpha\, d(x,x_0)\quad \text{ for all }\; x\in \G\,.
\end{equation}

\begin{lemma}\label{lemma41}
Let $p\geq 2$. Let $\xi$ be the function defined in \eqref{xi} and suppose \eqref{Vbound}. Let $\mathcal H:=\alpha^2-p\,V_0$, where $V_0$ is given in \eqref{Vbound} and assume
$$0<\alpha <\sqrt{p\,V_0}.$$
 Then, for every $e\in E$ and $x\in I_e$,
\begin{equation}\label{eq422}
(\xi_e'(x))^2-pV_e(x)\le\mathcal H,
\end{equation}
and
\begin{equation}\label{eq422b}
\xi''_e(x)+(\xi_e'(x))^2-pV_e(x)\le\mathcal H.
\end{equation}
\end{lemma}

\begin{proof}[Proof of Lemma \ref{lemma41}]
Inequalities \eqref{eq422} and \eqref{eq422b} are a direct consequence of the definition of $\xi$ in \eqref{xi} and assumption \eqref{Vbound}. In fact, for every $e\in E$
$$
(\xi_e'(x))^2= \alpha ^2,\quad \xi_e''(x)=0
$$
thus
$$
\xi''_e(x)+(\xi_e'(x))^2-pV_e(x)=(\xi_e'(x))^2-pV_e(x)\le \alpha ^2-p\,V_0\,.
$$
\end{proof}

Let $\hat \eta\in C^\infty([0, +\infty))$, $0\leq \hat \eta\leq 1,$
\[\hat\eta(t):=
\begin{cases}
1 & \text{ if } 0\leq t\leq 1,\\
0 & \text{ if } t>2\,.
\end{cases}
\]
We define 
\[\eta_R(x):=\hat\eta\left(\frac{d(x, x_0)}{R}\right) \quad \text{ for any } x\in \mathcal G\,.\]
Therefore $\eta_R\in C^2(\mathcal G)$, 
\begin{equation}\label{etadef}
\eta_R(x)=\begin{cases}1&\quad\text{if}\,\,\,x\in B_R(x_0)\\ 0&\quad\text{if}\,\,\,x\in (B_{2R}(x_0))^C.\end{cases}
\end{equation}
 Observe that, for some $C_1>0$, $C_2>0$,
 \begin{equation}\label{eta2}
 \begin{aligned}
 &|\eta_R'(x)|\le \frac {C_1}R\,\chi_{\{R\le d(x,x_0)\le 2R\}} \quad \text{for any}\quad x\in \G\,;\\
 &|\Delta\eta_R(x)|\le \frac {C_2}{R^2}\,\chi_{\{R\le d(x,x_0)\le 2R\}} \quad \text{for any}\quad x\in \G\,.
 \end{aligned}
 \end{equation}

\begin{proof}[Proof of Theorem \ref{teo1}]
Let $\xi$ be defined as in \eqref{xi} with $\alpha=2\beta$.

 \medskip

   We first prove item (i). Let $p\ge2$ and $0<\beta<\frac{\sqrt{p V_0}}{2}$. From Proposition \ref{prop41}, we get
\begin{equation}\label{eq13f}
\begin{aligned}
\int_{\G}&|u(x)|^p(\eta_R(x))^2e^{\xi(x)}\left[p\,V(x)-\left(\xi'(x)\right)^2\right]\,d\mu(x)\le4\int_{\G}|u(x)|^p\left(\eta_R'(x)\right)^2e^{\xi(x)}\,d\mu(x)\,.
\end{aligned}
\end{equation}
Since $\alpha=2\beta$, we have that $0<\alpha<\sqrt{pV_0}$ (see assumption \eqref{h1}), therefore we can apply Lemma \ref{lemma41} to the left-hand-side of \eqref{eq13f}, it yields
$$
\int_{\G}|u(x)|^p(\eta_R(x))^2e^{\xi(x)}\left[p\,V(x)-\left(\xi'(x)\right)^2\right]\,d\mu(x)\ge|\mathcal H|\int_{\G}|u(x)|^p(\eta_R(x))^2e^{\xi(x)}\,d\mu(x)\,.
$$
On the other hand, by the very definition of $\xi$ and estimate \eqref{eta2}, the right-hand-side of \eqref{eq13f} gives
$$
\begin{aligned}
\int_{\G}|u(x)|^p(\eta_R(x))^2e^{\xi(x)}&\left[p\,V(x)-\left(\xi'(x)\right)^2\right]\,d\mu(x)\le4\int_{\G}|u(x)|^p(\eta_R'(x))^2e^{-\alpha\, d(x,x_0)}\,d\mu(x)\\
&\le\frac{4C_1^2}{R^2}\int_{B_{2R}\setminus B_R}|u(x)|^pe^{-\alpha\, d(x,x_0)}\,d\mu(x)\\
&\le\frac{4C^2}{R^2}e^{-\alpha R}\int_{\G}|u(x)|^p\,d\mu(x)\,.
\end{aligned}
$$
Combining the above inequalities with \eqref{eq13f}, we get
$$
|\mathcal H|\int_{\G}|u(x)|^p(\eta(x))^2e^{\xi(x)}\,d\mu(x)\le\frac{4C^2}{R^2}e^{-\alpha R}\int_{B_{2R}}|u(x)|^p\,d\mu(x)\,.
$$
Now, since $\alpha=2\beta$ and due to assumption \eqref{growthcond1}, the latter reduces, for some $C>0$, to
\begin{equation}\label{eq426}
|\mathcal H|\int_{\G}|u(x)|^p(\eta(x))^2e^{\xi(x)}\,d\mu(x)\le \frac{4C_1^2C}{R^2}e^{-2\beta R+2\beta R}= \frac{4C_1^2C}{R^2}\,.
\end{equation}
We finally take the limit as $R\to+\infty$ in \eqref{eq426}, by dominated convergence theorem, since $\eta\to1$ and $|\mathcal H|>0$, we get
$$
\int_{\G}|u(x)|^pe^{\xi(x)}\,d\mu(x)\le0\,,
$$
which ensures that $u(x)=0$ for all $x\in \G$.
\medskip

We now verify item (ii). Let $1\le p<2$, $0<\beta<\frac{\sqrt{1+4pV_0}-1}{4}$ and let assumption \eqref{degAss} be verified.  We define $f(x)=\eta_R(x)\,e^{\xi(x)}$ for all $x\in \mathcal G$. Note that
 $$
 \Delta f(x)=\eta_R(x)\Delta\left(e^{\xi(x)}\right)+2\eta_R'(x)\xi'(x) e^{\xi(x)}+e^{\xi(x)}\Delta \eta_R(x)\,.
 $$
Then from Proposition \ref{prop42} we get
\begin{equation}\label{eq13g}
\begin{aligned}
\int_{\G}|u(x)|^p&\left[-\Delta\left(e^{\xi(x)}\right)+pV(x)e^{\xi(x)}\right]  \eta_R(x)\,d\mu(x)\\
&-\int_{\G}|u(x)|^p\Delta \eta_R(x)e^{\xi(x)} \,d\mu(x)-2\int_{\G}|u(x)|^p\eta_R'(x)\xi'(x) e^{\xi(x)}\,d\mu(x)\\
&\le-\sum_{x\in V}|u(x)|^pe^{\xi(x)}[\mathcal K(\eta_R)](x)-\sum_{x\in V}|u(x)|^pe^{\xi(x)}\eta_R(x)[\mathcal K(\xi)](x)
\end{aligned}
\end{equation}
Observe that $$\Delta\left(e^{\xi_e(x)}\right)=\left[\xi''_e(x)+(\xi'_e(x))^2\right]e^{\xi_e(x)},$$
therefore, the first term in the left-hand-side of \eqref{eq13g} gives
$$
\begin{aligned}
\int_{\G}|u(x)|^p&\left[-\Delta\left(e^{\xi(x)}\right)+pV(x)e^{\xi(x)}\right] \eta_R(x)\,d\mu(x)\\
&=\sum_{e\in E}\int_{I_e}|u_e(x)|^p\left[\xi''_e(x)+(\xi'_e(x))^2+pV_e(x)\right]e^{\xi_e(x)}\eta_{R,e}(x)\,dx.\end{aligned}
$$
Since $\alpha=2\beta$, it follows that $\alpha<\frac{\sqrt{pV_0}}{2}$ (see assumption \eqref{h1}), therefore we can apply Lemma \ref{lemma41} to the latter; this yields
$$
\int_{\G}|u(x)|^p\left[-\Delta\left(e^{\xi(x)}\right)+pV(x)e^{\xi(x)}\right] \eta_R(x)\,d\mu(x)\ge|\mathcal H|\int_{\G}|u(x)|^p\eta_R(x)e^{\xi(x)}\,d\mu(x)\,.
$$
Hence, \eqref{eq13g} reduces to
\begin{equation}\label{ineq1}
\begin{aligned}
|\mathcal H|\int_{\G}|u(x)|^p\eta_R(x)e^{\xi(x)}\,d\mu(x)\le I_{1,R}+I_{2,R}+I_{3,R}+I,
\end{aligned}
\end{equation}
where
$$
\begin{aligned}
I_{1,R}&=\int_{\G}|u(x)|^p\Delta \eta_R(x)e^{\xi(x)} \,d\mu(x)\\
I_{2,R}&=2\int_{\G}|u(x)|^p\eta_R^{(1)}(x)\xi^{(1)}(x) e^{\xi(x)}\,d\mu(x)\\
I_{3,R}&=-\sum_{x\in V}|u(x)|^pe^{\xi(x)}[\mathcal K(\eta_R)](x)\\
I&=-\sum_{x\in V}|u(x)|^pe^{\xi(x)}\eta_R(x)[\mathcal K(\xi)](x)\,.
\end{aligned}
$$
The thesis follows if we show that $I_{i,R}\to0$ as $R\to\infty$ for $i=1,2,3$ and $I\le0$. Since $\alpha=2\beta$, due to assumption \eqref{growthcond1} and \eqref{eta2}, we have, for some $C>0$,
\begin{equation}\label{limit}
\begin{aligned}
|I_{1,R}|&\le\int_{\G}|u(x)|^p|\Delta \eta_R(x)|e^{\xi(x)} \,d\mu(x) \le\frac{C_2}{R^2}\int_{B_{2R}\setminus B_R}|u(x)|^p e^{-\alpha d(x,x_0)} \,d\mu(x)\\
&\le \frac{C_2}{R^2} e^{-\alpha R}\int_{B_{2R}}|u(x)|^p \,d\mu(x)\le \frac{C\,C_2}{R^2}e^{-\alpha R+2\beta R}=\frac{C\,C_2}{R^2}\,;\\
|I_{2,R}|&\le 2\int_{\G}|u(x)|^p|\eta_R'(x)||\xi'(x)| e^{\xi(x)}\,d\mu(x)\le \frac{2 \beta C_1}{R}\int_{\G}|u(x)|^pe^{\xi(x)}\,d\mu(x)\\
&\le \frac{2 \beta C_1}{R} e^{-\alpha R}\int_{B_{2R}}|u(x)|^p \,d\mu(x)\le \frac{2 \beta C_1}{R}e^{-\alpha R+2\beta R}= \frac{2\beta C\, C_1}{R}\,.\\
\end{aligned}
\end{equation}
Moreover, from the very definition of $\mathcal K$ in \eqref{K} and $\eta_R$,
\begin{equation}\label{limit2}
|I_{3,R}|\le\sum_{x\in V}|u(x)|^pe^{\xi(x)}|[\mathcal K(\eta_R)](x)|\le \frac{ C_1}{R}\sum_{x\in V}|u(x)|^pe^{\xi(x)}\,.
\end{equation}
Finally, we observe that, since $\xi'\le0$, due to \eqref{e502u} and to assumption \eqref{degAss}, we obtain
\begin{equation*}
[\mathcal K(\xi)](x)=\sum_{k=1}^{\operatorname{deg}_x^+}\xi'(x)-\sum_{k=1}^{\operatorname{deg}_x^-}\xi'(x)=\xi'(x)\left(\operatorname{deg}_x^+-\operatorname{deg}_x^-\right)\ge0\quad\text{for every $x\in \mathcal V$},
\end{equation*}
therefore
\begin{equation}\label{boundary2}
I=-\sum_{x\in V}|u(x)|^pe^{\xi(x)}\eta_R(x)[\mathcal K(\xi)](x)\le0\,.
\end{equation}
By combining \eqref{ineq1}, \eqref{limit}, \eqref{limit2} and \eqref{boundary2}, we get
$$
|\mathcal H|\int_{\G}|u(x)|^p\eta_R(x)e^{\xi(x)}\,d\mu(x)\le\frac{C\,C_2}{R^2}+\frac{2 \beta C\, C_1}{R}+\frac{ C_1}{R}\sum_{x\in V}|u(x)|^pe^{\xi(x)}\,.
$$
We finally let $R\to+\infty$, we obtain
$$
|\mathcal H|\int_{\G}|u(x)|^pe^{\xi(x)}\,d\mu(x)\le0\,.
$$
Since $|\mathcal H|>0$, the thesis follows.
\end{proof}

\subsection{Vanishing potential}

For any $\sigma>0$, we define the function
 \begin{equation}\label{zeta}
\zeta(x):=\left[ d(x,x_0)+k\right]^{-\sigma}\quad \text{ for all }\; x\in \G\,.
\end{equation}

\begin{lemma}\label{lemma71}
Let $p\geq 1$. Let $\zeta$ be the function defined in \eqref{zeta} and suppose \eqref{Vbound2}. Then, for every $e\in E$ and $x\in e$,
\begin{itemize}
\item if
$$
0<\sigma<\sqrt{pV_0}\,,
$$
then there exists $K>0$ such that
\begin{equation}\label{eq71}
\left[\frac{\zeta_e'(x)}{\zeta_e(x)}\right]^2-pV_e(x)\le -K\,[d(x,x_0)+k]^{-2};
\end{equation}
\item if
$$
0<\sigma<\frac{\sqrt{1+4pV_0}-1}{2},
$$
then there exists $H>0$ such that
\begin{equation}\label{eq72}
\zeta''_e(x)-p V_e(x)\le -H [d(x,x_0)+k]^{-2}.
\end{equation}
\end{itemize}
Here $V_0$ is the same as in \eqref{Vbound2}.
\end{lemma}

\begin{proof}[Proof of Lemma \ref{lemma71}]
We compute, for every $e\in E$, $x\in I_e$,
$$
\zeta_e'(x)=-\sigma [d(x,x_0)+k]^{-\sigma-1}.
$$
Therefore, by \eqref{Vbound2}
$$
\begin{aligned}
\left[\frac{\zeta_e'(x)}{\zeta_e(x)}\right]^2-pV_e(x)&=\sigma^2 \frac{[d(x,x_0)+k]^{-2\sigma-2}}{[d(x,x_0)+k]^{-2\sigma}}-pV_e(x)\\
&\le\sigma^2 [d(x,x_0)+k]^{-2}-p\,V_0[d(x,x_0)+k]^{-\theta}\\
&\le-(p\,V_0-\sigma^2)[d(x,x_0)+k]^{-2}\,.
\end{aligned}
$$
Inequality \eqref{eq71} follows due to assumption $p\,V_0-\sigma^2>0$.
Furthermore, due to \eqref{Vbound2} with $0<\theta\le2$
$$
\begin{aligned}
\zeta''_e(x)-pV_e(x)&\le\sigma(\sigma+1)[d(x,x_0)+k]^{-\sigma-2}-pV_e(x)\\
&\le\sigma(\sigma+1)[d(x,x_0)+k]^{-\sigma-2}-p\,V_0[d(x,x_0)+k]^{-\theta}\\
&\le-(p\,V_0-\sigma^2-\sigma)[d(x,x_0)+k]^{-2}\,.
\end{aligned}
$$
Inequality \eqref{eq72} follows due to assumption $p\,V_0-2\sigma^2-\sigma>0$.
\end{proof}

\begin{proof}[Proof of Theorem \ref{teo3}]
Let $\zeta$ be defined as in \eqref{zeta} with $\sigma=\lambda$.
Let $\eta_R\in\mathfrak F$ be defined as in \eqref{etadef}.
  We first prove item (i). Let $p\ge2$ and $\lambda<\sqrt{pV_0}$. From Corollary \ref{cor41}, we get
\begin{equation}\label{eq725}
\begin{aligned}
\int_{\G}&|u(x)|^p\eta^2(x){\zeta(x)}\left[p\,V(x)-\left(\frac{\zeta'(x)}{\zeta(x)}\right)^2\right]\,d\mu(x)\le4\int_{\G}|u(x)|^p\left(\eta'(x)\right)^2{\zeta(x)}\,d\mu(x)\\
\end{aligned}
\end{equation}
Since $\sigma=\lambda$, we have $\sigma<\sqrt{pV_0}$ (see assumption \eqref{h3}), therefore we can apply Lemma \ref{lemma71} to the left-hand-side of \eqref{eq725}, it yields
$$
\int_{\G}|u(x)|^p\eta^2(x){\zeta(x)}\left[p\,V(x)-\left(\frac{\zeta'(x)}{\zeta(x)}\right)^2\right]\,d\mu(x)\ge K\int_{\G}|u(x)|^p(\eta(x))^2{\zeta(x)}(d(x,x_0)+k)^{-2}\,d\mu(x)\,.
$$
On the other hand, the right-hand-side of \eqref{eq725}, using the very definition of $\zeta$ and estimate \eqref{etadef}, gives
$$
\begin{aligned}
\int_{\G}|u(x)|^p\eta^2(x)&{\zeta(x)}\left[p\,V(x)-\left(\frac{\zeta'(x)}{\zeta(x)}\right)^2\right]\,d\mu(x)\le4\int_{\G}|u(x)|^p(\eta'(x))^2(d(x,x_0)+k)^{-\sigma}\,d\mu(x)\\
&\le\frac{4C_1^2}{R^2}\int_{B_{2R}\setminus B_R}|u(x)|^p(d(x,x_0)+k)^{-\sigma}\,d\mu(x)\\
&\le\frac{4C^2}{R^2}(R+k)^{-\sigma}\int_{B_{2R}}|u(x)|^p\,d\mu(x)\,.
\end{aligned}
$$
Thus, combining the above estimates with \eqref{eq725}, we get
$$
K\int_{\G}|u(x)|^p(\eta(x))^2{\zeta(x)}(d(x,x_0)+k)^{-2}\,d\mu(x)\le\frac{4C^2}{R^2}(R+k)^{-\sigma}\int_{B_{2R}}|u(x)|^p\,d\mu(x)\,.
$$
Now, since $\sigma=\lambda$ and due to assumption \eqref{growthcond2},
the latter reduces, for some $\tilde C>0$, to
\begin{equation}\label{eq726}
K\int_{\G}|u(x)|^p(\eta(x))^2{\zeta(x)}(d(x,x_0)+k)^{-2}\,d\mu(x)\le\frac{4C_1^2\tilde C}{R^2}\left(\frac{2R+k}{R+k}\right)^{\lambda}\le\frac{8C_1^2\tilde C}{R^2}\,.
\end{equation}
We take the limit as $R\to+\infty$ in \eqref{eq726}, by dominated convergence theorem, since $\eta\to1$ and $K>0$, we get
$$
\int_{\G}|u(x)|^p{\zeta(x)}(d(x,x_0)+k)^{-2}\,d\mu(x)\le0\,,
$$
which ensures that $u(x)=0$ for all $x\in \G$ because $\zeta>0$ and $(d(\cdot,x_0)+k)^{-2}>0$ in $\mathcal G$.
\medskip

We now verify item (ii). Let $1\le p<2$, $0<\lambda<\frac{\sqrt{1+4pV_0}-1}{2}$ and let assumption \eqref{degAss} be verified.  We define $g(x)=\eta_R(x)\,{\zeta(x)}$, for all $x\in \mathcal G$. Note that
 $$
 \Delta g(x)=\eta_R(x)\Delta\zeta(x)+2\eta_R'(x)\zeta'(x) +{\zeta(x)}\Delta \eta_R(x)\,.
 $$
Then from Proposition \ref{prop42} we get
\begin{equation}\label{eq13h}
\begin{aligned}
\int_{\G}|u(x)|^p&\left[-\Delta\zeta(x)+pV(x)\zeta(x)\right]  \eta_R(x)\,d\mu(x)\\
&-\int_{\G}|u(x)|^p\Delta \eta_R(x)\zeta(x)\,d\mu(x)-2\int_{\G}|u(x)|^p\eta_R'(x)\zeta'(x)\,d\mu(x)\\
&\le-\sum_{x\in V}|u(x)|^p\zeta(x)[\mathcal K(\eta_R)](x)-\sum_{x\in V}|u(x)|^p\zeta(x)\eta_R(x)[\mathcal K(\xi)](x)
\end{aligned}
\end{equation}
Observe that $\Delta\zeta_e(x)=\zeta''_e(x),$
therefore \eqref{eq13h} reads
$$
\begin{aligned}
\int_{\G}|u(x)|^p&\left[-\zeta''(x)+pV(x)\zeta(x)\right]  \eta_R(x)\,d\mu(x)\\
&-\int_{\G}|u(x)|^p\Delta \eta_R(x)\zeta(x)\,d\mu(x)-2\int_{\G}|u(x)|^p\eta_R'(x)\zeta'(x)\,d\mu(x)\\
&\le-\sum_{x\in V}|u(x)|^p\zeta(x)[\mathcal K(\eta_R)](x)-\sum_{x\in V}|u(x)|^p\eta_R(x)[\mathcal K(\zeta)](x)
\end{aligned}
$$
Since $\sigma=\lambda$ we also have $0<\lambda<\frac{\sqrt{1+4pV_0}-1}{2}$ (see assumption \eqref{h4} of the theorem), therefore we can apply Lemma \ref{lemma71} to the latter, this yields
\begin{equation}\label{ineq2}
\begin{aligned}
K\int_{\G}|u(x)|^p\eta_R(x)[d(x,x_0)+k]^{-2}\,d\mu(x)\le J_{1,R}+J_{2,R}+J_{3,R}+J,
\end{aligned}
\end{equation}
where
$$
\begin{aligned}
J_{1,R}&=\int_{\G}|u(x)|^p\Delta \eta_R(x)\zeta(x)\,d\mu(x)\\
J_{2,R}&=2\int_{\G}|u(x)|^p\eta_R'(x)\zeta'(x)\,d\mu(x)\\
J_{3,R}&=-\sum_{x\in V}|u(x)|^p\zeta(x)[\mathcal K(\eta_R)](x)\\
J&=-\sum_{x\in V}|u(x)|^p\eta_R(x)[\mathcal K(\zeta)](x)\,.
\end{aligned}
$$
The thesis follows if we show that $J_{i,R}\to0$ as $R\to\infty$ for $i=1,2,3$ and $J\le0$. Due to assumptions \eqref{growthcond2} and \eqref{etadef}, since $\sigma=\lambda$, we have, for some $C>0$,
\begin{equation}\label{limitsecond}
\begin{aligned}
|J_{1,R}|&\le\int_{\G}|u(x)|^p|\Delta \eta_R(x)|\zeta(x)\,d\mu(x) \le\frac{C_2}{R^2}\int_{B_{2R}\setminus B_R}|u(x)|^p [d(x,x_0)+k]^{-\sigma} \,d\mu(x)\\
&\le\frac{C_2}{R^2} (R+k)^{-\sigma}\int_{B_{2R}}|u(x)|^p \,d\mu(x) \le\frac{C_2}{R^2} \left(\frac{2R+k}{R+k}\right)^{\lambda}\le\frac{2C\,C_2}{R^2}\,;\\
|J_{2,R}|&\le 2\int_{\G}|u(x)|^p|\eta_R'(x)||\zeta'(x)|\,d\mu(x)\\
&\le \frac{2 \lambda C_1}{R}\int_{B_{2R}\setminus B_R}|u(x)|^p[d(x,x_0)+k]^{-\sigma-1}\,d\mu(x)\\
&\le \frac{2 \lambda C_1}{R}(R+k)^{-\sigma-1}\int_{B_{2R}}|u(x)|^p\,d\mu(x)\le\frac{2 \lambda C_1}{R(R+k)}\left(\frac{2R+k}{R+k}\right)^{\lambda}\le \frac{4\lambda C\, C_1}{R^2}\\
\end{aligned}
\end{equation}
Moreover, from the very definition of $\mathcal K$ in \eqref{K} and $\eta_R$,
\begin{equation}\label{limitsecond2}
|J_{3,R}|\le\sum_{x\in V}|u(x)|^p\zeta(x)|[\mathcal K(\eta_R)](x)|\le \frac{ C_1}{R}\sum_{x\in V}|u(x)|^p\zeta(x)\,.
\end{equation}
Finally, we observe that, since $\zeta'\le0$, due to \eqref{e502u} and assumption \eqref{degAss}
\begin{equation*}
[\mathcal K(\zeta)](x)=\sum_{k=1}^{\operatorname{deg}_x^+}\zeta'(x)-\sum_{k=1}^{\operatorname{deg}_x^-}\zeta'(x)=\zeta'(x)\left(\operatorname{deg}_x^+-\operatorname{deg}_x^-\right)\ge0\quad\text{for every $x\in V$},
\end{equation*}
therefore
\begin{equation}\label{boundarysecond}
J=-\sum_{x\in V}|u(x)|^p\eta_R(x)[\mathcal K(\zeta)](x)\le0\,.
\end{equation}
By combining \eqref{ineq2}, \eqref{limitsecond}, \eqref{limitsecond2} and \eqref{boundarysecond}, we get
$$
K\int_{\G}|u(x)|^p\eta_R(x)[d(x,x_0)+k]^{-2}\,d\mu(x)\le\frac{2C\,C_2}{R^2}+\frac{4\lambda C\, C_1}{R^2}+\frac{ C_1}{R}\sum_{x\in V}|u(x)|^p\zeta(x)\,.
$$
We finally let $R\to+\infty$: by dominated convergence theorem, since $\eta_R\to1$, we obtain
$$
K\int_{\G}|u(x)|^p[d(x,x_0)+k]^{-2}\,d\mu(x)\le0\,.
$$
Since $K>0$ and $[d(x,x_0)+k]^{-2}>0$, the thesis follows.

\end{proof}

\section{Proofs of the main results for the parabolic problem}\label{sec6}\setcounter{equation}{0}

This section is devoted to the proof of Theorems \ref{teo5} and \ref{teo7}. Therefore, we need to consider separately the cases when the density $\rho$ satisfies either assumption \eqref{rhobound} or \eqref{rhobound2}.

\subsection{Density bounded away from zero}

 We first define the test function that we will use to prove our main results. This test function is selected depending on the class of uniqueness which we are considering, which in turn depends on the assumption made on $\rho$, i.e. \eqref{rhobound}. Define
\begin{equation}\label{omega}
\omega(x,t):=-\gamma t+\xi(x)\quad \text{ for all }\; x\in \G,\,\, t\in[0,T]\,,
\end{equation}
for any $\gamma\ge0$ and $\xi$ as in \eqref{xi}, i.e. $\xi(x)=-\alpha d(x,x_0)$ with $\alpha\ge0$. The following lemma which will play a key role in the proof of Theorem \ref{teo5}.

\begin{lemma}\label{lemma61}
Let $p\geq 2$. Let $\omega$ be the function defined in \eqref{omega}. Suppose that \eqref{rhobound} holds, and that
$$\gamma\ge\frac{\alpha ^2}{\rho_0},$$
with $\rho_0$ as in \eqref{rhobound}. Then, for every $e\in E$ and $x\in I_e$,
\begin{align}
&\rho(x)\partial_t\omega_e+(\omega_e'(x,t))^2\le0,\label{eq61}\\
&\rho(x)\partial_t\omega_e+\omega''_e(x,t)+(\omega_e'(x,t))^2\le0.\label{eq62}
\end{align}
\end{lemma}

\begin{proof}[Proof of Lemma \ref{lemma61}]
Inequalities \eqref{eq61} and \eqref{eq62} are a direct consequence of the definition of $\omega$ in \eqref{omega}, $\xi$ in \eqref{xi} and assumption \eqref{rhobound}. In fact, for every $e\in E$, $x\in E$ and $t\in [0,T]$,
$$
(\omega_e'(x,t))^2= \alpha ^2,\quad \omega_e''(x,t)=0, \quad \partial_t\omega(x,t)=-\gamma.
$$
Thus, by the assumption on $\gamma$ and by \eqref{rhobound},
$$
\begin{aligned}
\rho(x)&\partial_t\omega_e+\omega''_e(x,t)+(\omega_e'(x,t))^2\\
&=\rho(x)\partial_t\omega_e+(\omega_e'(x,t))^2= -\gamma\rho(x)+\alpha^2\le -\gamma\rho_0+\alpha^2\le0\,.
\end{aligned}
$$
\end{proof}

\begin{proof}[Proof of Theorem \ref{teo5}]
Let $\omega$ be defined as in \eqref{omega} with $\alpha=2\beta$ and $\gamma\ge\frac{\alpha^2}{\rho_0}$. Moreover, let $\eta_R$ be defined as in \eqref{etadef}.
 
 \medskip

 We first prove item (i). Let $p\ge2$. From Proposition \ref{prop51} with the latter choices of $\omega$ and $\eta_R$, we get, for some $\tau\in(0,T]$
\begin{equation}\label{eq65}
\begin{aligned}
\int_{\G}\,\rho(x)\,&|u(x,\tau)|^pe^{\omega(x,\tau)}\eta_R^2(x)\,d\mu(x)\le 4\int_0^\tau\int_{\G}|u(x,t)|^p(\eta_R'(x))^2e^{\omega(x,t)}\, d\mu(x)dt\\
&\quad+\int_0^\tau \int_{\G}|u(x,t)|^p\left[\rho(x)\,\partial_t \omega(x,t)+(\omega'(x,t))^2\right]e^{\omega(x,t)}\eta_R^2(x)\,d\mu(x)dt\,.
\end{aligned}
\end{equation}
Since $\gamma\ge\frac{\alpha^2}{\rho_0}$, we can apply Lemma \ref{lemma61} to the last term of \eqref{eq65}, it yields
$$
\int_{\G}\,\rho(x)\,|u(x,\tau)|^pe^{\omega(x,\tau)}\eta_R^2(x)\,d\mu(x)\le 4\int_0^\tau\int_{\G}|u(x,t)|^p(\eta_R'(x))^2e^{\omega(x,t)}\, d\mu(x)dt\,.
$$
Then, by the very definition of $\omega$ and estimate \eqref{eta2}, the latter gives
$$
\begin{aligned}
\int_{\G}\,\rho(x)\,|u(x,\tau)|^p&e^{\omega(x,\tau)}\eta_R^2(x)\,d\mu(x)\le
4\int_0^\tau\int_{\G}|u(x,t)|^p(\eta_R'(x))^2e^{-\gamma t-\alpha d(x,x_0)}\, d\mu(x)\,dt\\
&\le\frac{4C_1^2}{R^2}\int_0^\tau\int_{B_{2R}\setminus B_R}|u(x,t)|^pe^{-\alpha\, d(x,x_0)}\,d\mu(x)\,dt\\
&\le\frac{4C_1^2}{R^2}e^{-\alpha\, R}\int_0^\tau\int_{B_{2R}}|u(x,t)|^p\,d\mu(x)\,dt\,.
\end{aligned}
$$
Since $\alpha=2\beta$ and due to assumption \eqref{growthcond3},
the latter reduces, for some $C>0$, to
\begin{equation}\label{eq66}
\int_{\G}\,\rho(x)\,|u(x,\tau)|^pe^{\omega(x,\tau)}\eta_R^2(x)\,d\mu(x)\le \frac{4C_1^2C}{R^2}e^{-2\beta R+2\beta R}= \frac{4C_1^2C}{R^2}\,.
\end{equation}
We finally take the limit as $R\to+\infty$ in \eqref{eq66}, by dominated convergence theorem, since $\eta_R\to1$, we get
$$
\int_{\G}\,\rho(x)\,|u(x,\tau)|^pe^{\omega(x,\tau)}\,d\mu(x)\le0\,.
$$
Observe that $\rho>0$ (see assumption \eqref{rhobound}), $e^\omega>0$ in $\G\times(0,T)$, therefore it necessarily holds that $u(x,\tau)=0$ for all $x\in \G$. By arbitrariness of $\tau\in(0,T]$, we get the thesis.
%
%
\medskip

We now verify item (ii). Let $1\le p<2$ and let assumption \eqref{degAss} be verified.  Let $\tau\in(0,T]$, we define  $v(x,t)=\eta_R(x)\,e^{\omega(x,t)}$ ($x\in \mathcal G$, $t\in(0,\tau]$). 
Note that
 $$
 \Delta v(x,t)=\eta_R(x)\Delta\left(e^{\omega(x,t)}\right)+2\eta_R'(x)\omega'(x,t) e^{\omega(x,t)}+e^{\omega(x,t)}\Delta \eta_R(x)\,.
 $$
We can apply Proposition \ref{prop52} with this choice of $v$, thus we get
\begin{equation}\label{eq67}
\begin{aligned}
\int_0^\tau&\int_{\G}|u(x,t)|^p\left[\rho(x)\partial_t\left(e^{\omega(x,t)}\right)+\Delta \left(e^{\omega(x,t)}\right)\right] \eta_R(x)\,d\mu(x)dt\\
&+\int_0^\tau\int_{\G}|u(x,t)|^p\Delta \left(\eta_R(x)\right) e^{\omega(x,t)}\,d\mu(x)\,dt+2\int_0^\tau\int_{\G}|u(x,t)|^p\eta'_R(x) \omega'(x,t)e^{\omega(x,t)}\,d\mu(x)dt\\
&\ge\int_{\G}\rho(x)|u(x,\tau)|^p\eta_R(x)e^{\omega(x,\tau)}\,d\mu(x)+\int_0^\tau\sum_{x\in V}|u(x,t)|^p[\mathcal K(\eta_Re^{\omega})](x,t)dt\,.\\
\end{aligned}
\end{equation}
Observe that $$\Delta\left(e^{\omega_e(x,t)}\right)=\left[\omega''_e(x,t)+(\omega'_e(x,t))^2\right]e^{\omega_e(x,t)},$$ therefore, the first term in the left-hand-side of \eqref{eq67} gives
$$
\begin{aligned}
\int_0^\tau&\int_{\G}|u(x,t)|^p\left[\rho(x)\partial_t\left(e^{\omega(x,t)}\right)+\Delta \left(e^{\omega(x,t)}\right)\right] \eta_R(x)\,d\mu(x)dt\\
&=\sum_{e\in E}\int_{I_e}|u_e(x,t)|^p\left[\partial_t\omega_e(x,t)+\omega''_e(x,t)+(\omega'_e(x,t))^2\right]e^{\omega_e(x,t)}\eta_{R,e}(x)\,dx.\end{aligned}
$$
Since $\gamma\ge\frac{\alpha^2}{\rho_0}$, we can apply Lemma \ref{lemma61} to the latter, this yields
$$
\int_0^\tau\int_{\G}|u(x,t)|^p\left[\rho(x)\partial_t\left(e^{\omega(x,t)}\right)+\Delta \left(e^{\omega(x,t)}\right)\right] \eta_R(x)\,d\mu(x)dt\le0
$$
Hence, \eqref{eq67} reduces to
\begin{equation}\label{eq68}
\begin{aligned}
\int_0^\tau\int_{\G}&|u(x,t)|^p\Delta \left(\eta_R(x)\right) e^{\omega(x,t)}\,d\mu(x)dt\\
&\quad+2\int_0^\tau\int_{\G}|u(x,t)|^p\eta'_R(x) \omega'(x,t)e^{\omega(x,t)}\,d\mu(x)dt\\
&\ge\int_{\G}\rho(x)|u(x,\tau)|^p\eta_R(x)e^{\omega(x,\tau)}\,d\mu(x)+\int_0^\tau\sum_{x\in V}|u(x,t)|^p[\mathcal K(\eta_Re^{\omega})](x,t)\,dt\,.\\
&\ge\int_{\G}\rho(x)|u(x,\tau)|^p\eta_R(x)e^{\omega(x,\tau)}\,d\mu(x)\\
&\quad+\int_0^\tau\sum_{x\in V}|u(x,t)|^p\eta_R(x)e^{\omega}(x,t)[\mathcal K({\omega})](x,t)\,dt\\
&\quad+\int_0^\tau\sum_{x\in V}|u(x,t)|^pe^{\omega}(x,t)[\mathcal K(\eta_R)](x)\,dt\,.\\
\end{aligned}
\end{equation}
Defining
$$
\begin{aligned}
I_{1,R}&=\int_0^\tau\int_{\G}|u(x,t)|^p\Delta \left(\eta_R(x)\right) e^{\omega(x,t)}\,d\mu(x)dt\\
I_{2,R}&=2\int_0^\tau\int_{\G}|u(x,t)|^p\eta'_R(x) \omega'(x,t)e^{\omega(x,t)}\,d\mu(x)dt\\
I_{3,R}&=-\int_0^\tau\sum_{x\in V}|u(x,t)|^pe^{\omega(x,t)}[\mathcal K(\eta_R)](x)\,dt\,\\
I&=-\int_0^\tau\sum_{x\in V}|u(x,t)|^p\eta_R(x)e^{\omega(x,t)}[\mathcal K({\omega})](x,t)\,dt\,;
\end{aligned}
$$
we can rearrange \eqref{eq68} as
\begin{equation}\label{eq68a}
\int_{\G}\rho(x)|u(x,\tau)|^p\eta_R(x)e^{\omega(x,\tau)}\,d\mu(x)\le I_{1,R}+I_{2,R}+I_{3,R}+I\,.
\end{equation}
The thesis follows if we show that $I_{i,R}\to0$ as $R\to\infty$ for $i=1,2,3$ and $I\le0$. Due to assumptions \eqref{growthcond3} and \eqref{eta2}, and since $\alpha=2\beta$, we have, for some $C>0$,
\begin{equation}\label{eq69}
\begin{aligned}
|I_{1,R}|&\le\int_0^\tau\int_{\G}|u(x,t)|^p|\Delta \eta_R(x)|e^{-\gamma t-\alpha d(x,x_0)} \,d\mu(x)dt \\
&\le\frac{C_2}{R^2}\int_0^\tau\int_{B_{2R}\setminus B_R}|u(x,t)|^p e^{-\alpha d(x,x_0)} \,d\mu(x)dt\\
&\le\frac{C_2}{R^2}e^{-\alpha R} \int_0^\tau\int_{B_{2R}}|u(x,t)|^p \,d\mu(x)dt\le\frac{CC_2}{R^2}e^{-\alpha R+2\beta R} =\frac{CC_2}{R^2}\,;\\
|I_{2,R}|&\le 2\int_0^\tau\int_{\G}|u(x,t)|^p|\eta_R'(x)|| \omega'(x,t)| e^{-\gamma t-\alpha d(x,x_0)}\,d\mu(x)dt\\
&\le \frac{2 \beta C_1}{R}\int_0^\tau\int_{B_{2R}\setminus B_R}|u(x,t)|^pe^{-\alpha d(x,x_0)}\,d\mu(x)dt\\
&\le \frac{2 \beta C_1}{R}e^{-\alpha R} \int_0^\tau\int_{B_{2R}}|u(x,t)|^p \,d\mu(x)dt\le \frac{2\beta C C_1}{R}e^{-\alpha R+2\beta R} =\frac{2\beta C C_1}{R}\\
\end{aligned}
\end{equation}
Moreover, from the very definition of $\mathcal K$ in \eqref{K} and $\eta_R$,
\begin{equation}\label{eq610}
|I_{3,R}|\le\int_0^\tau\sum_{x\in V}|u(x,t)|^pe^{\omega(x,t)}|[\mathcal K(\eta_R)](x)|\,dt\le \frac{ C_1}{R}\int_0^\tau\sum_{x\in V}|u(x,t)|^pe^{\omega(x,t)}\,dt\,.
\end{equation}
Finally, we observe that, since $\omega'\le0$, due to \eqref{e502u} and assumption \eqref{degAss}, for every $t\in(0,\tau)$,
\begin{equation*}
[\mathcal K(\omega)](x,t)=\sum_{k=1}^{\operatorname{deg}_x^+}\omega'(x,t)-\sum_{k=1}^{\operatorname{deg}_x^-}\omega'(x,t)=\omega'(x,t)\left(\operatorname{deg}_x^+-\operatorname{deg}_x^-\right)\ge0\quad\text{for every}\,\, x\in V,
\end{equation*}
therefore
\begin{equation}\label{eq611}
I=-\int_0^\tau\sum_{x\in V}|u(x,t)|^pe^{\omega(x,t)}\eta_R(x)[\mathcal K(\omega)](x,t)dt\le0\,.
\end{equation}
By combining \eqref{eq69}, \eqref{eq610}, \eqref{eq611} and \eqref{eq68a}, we get
$$
\int_{\G}\rho(x)|u(x,\tau)|^p\eta_R(x)e^{\omega(x,\tau)}\,d\mu(x)\le\frac{C\,C_2}{R^2}+\frac{2 \beta C\, C_1}{R}+\frac{ C_1}{R}\int_0^\tau\sum_{x\in V}|u(x,t)|^pe^{\omega(x,t)}\,dt\,.
$$
We finally let $R\to+\infty$, we obtain
$$
\int_{\G}\rho(x)|u(x,\tau)|^pe^{\omega(x,\tau)}\,d\mu(x)\le0\,.
$$
Since $\rho>0$ (see assumption \eqref{rhobound}) and by arbitrariness of $\tau$, the thesis follows.
\end{proof}

\subsection{Vanishing density}
In this second part of our proofs we deal with the following test function
 \begin{equation}\label{kappa}
\kappa(x,t):=e^{-\gamma t}\zeta(x)\quad \text{ for all }\; x\in \G,\,\,t\in(0,T]\,.
\end{equation}
where $\gamma>0$ and $\zeta$ is given in \eqref{zeta}, i.e. $\zeta(x)=\left[ d(x,x_0)+k\right]^{-\sigma}$ for some $\sigma>0$ and $k\ge1$. Similarly to the previous subsection, also in this case it is useful to show that $\kappa$ satisfies some inequalities. For this reason we state the following Lemma which will play a key role in the proof of Theorem \ref{teo7}.

\begin{lemma}\label{lemma62}
Let $p\geq 1$. Let $\kappa $ be the function defined in \eqref{kappa}. Suppose \eqref{rhobound2} and that
$$\gamma\ge\frac{\sigma(\sigma+1)}{\rho_0},$$
for $\rho_0$ as in \eqref{rhobound2}. Then, for every $e\in E$, $x\in I_e$ and $t\in(0,T]$,
\begin{align}
&\rho(x)\partial_t\kappa_e(x,t)+\frac{(\kappa_e'(x,t))^2}{\kappa_e(x,t)}\le 0,\label{eq613}\\
& \rho(x)\partial_t\kappa_e(x,t)+\kappa''_e(x,t)\le 0.\label{eq614}
\end{align}
\end{lemma}

\begin{proof}[Proof of Lemma \ref{lemma62}]
We compute, for every $e\in E$ and $t\in(0,T]$,
$$
\begin{aligned}
&\partial_t\kappa_e(x,t)=-\gamma e^{-\gamma t}\zeta_e(x)=-\gamma\kappa_e(x,t);\\
&\kappa_e'(x,t)=e^{-\gamma t} \zeta'(x)=-\sigma e^{-\gamma t}\left[ d(x,x_0)+k\right]^{-\sigma-1}\operatorname{sgn}(x-x_0)\\
&\quad\quad \quad =-\sigma \left[ d(x,x_0)+k\right]^{-1}\operatorname{sgn}(x-x_0)\kappa_e(x,t);\\
&\kappa_e''(x,t)=e^{-\gamma t} \zeta''(x)=\sigma(\sigma+1) e^{-\gamma t}\left[ d(x,x_0)+k\right]^{-\sigma-2}\\
&\quad\quad \quad =\sigma(\sigma+1)\left[ d(x,x_0)+k\right]^{-2}\kappa_e(x,t).\\
\end{aligned}
$$
Therefore, by \eqref{rhobound2} with $0<\theta\le2$,
$$
\begin{aligned}
\rho(x)\partial_t\kappa_e(x,t)+\frac{(\kappa_e'(x,t))^2}{\kappa_e(x,t)}&=-\gamma\rho_e(x)\kappa_e(x,t)+\sigma^2 \left[ d(x,x_0)+k\right]^{-2}\kappa_e(x,t)\\
&\le\kappa_e(x,t)\left\{-\gamma\rho_0\left[ d(x,x_0)+k\right]^{-\theta}+\sigma^2\left[ d(x,x_0)+k\right]^{-2}\right\}\,.\\
\end{aligned}
$$
Inequality \eqref{eq613} follows because, by assumption, $\gamma\ge\frac{\sigma(\sigma+1)}{\rho_0}\left(\ge\frac{\sigma^2}{\rho_0}\right)$.

\noindent Furthermore, due to \eqref{rhobound2} with $0<\theta\le2$,
$$
\begin{aligned}
 \rho(x)\partial_t\kappa_e(x,t)+\kappa''_e(x)&=-\gamma\rho_e(x)\kappa_e(x,t)+\sigma(\sigma+1)[d(x,x_0)+k]^{-2}\kappa_e(x,t)\\
&\le\kappa_e(x,t)\left\{-\gamma\rho_0\left[ d(x,x_0)+k\right]^{-\theta}+\sigma(\sigma+1)[d(x,x_0)+k]^{-2}\right\}\,.\\
\end{aligned}
$$
Inequality \eqref{eq614} follows because, by assumption, $\gamma\ge\frac{\sigma(\sigma+1)}{\rho_0}$.
\end{proof}

\begin{proof}[Proof of Theorem \ref{teo7}]
Let $\kappa$ be defined as in \eqref{kappa} with $\sigma=\lambda$ and $\gamma\ge\frac{\sigma(\sigma+1)}{\rho_0}$. Let $\eta_R$ be defined as in \eqref{etadef}.

 We first prove item (i). Let $p\ge2$. From Corollary \ref{cor51} with the latter choices of $\kappa$ and $\eta_R$, we get, for some $\tau\in(0,T]$
\begin{equation}\label{eq616}
\begin{aligned}
\int_{\G}\,\rho(x)\,&|u(x,\tau)|^p{\kappa(x,\tau)}\eta_R^2(x)\,d\mu(x)\le 4\int_0^\tau\int_{\G}|u(x,t)|^p(\eta_R'(x))^2{\kappa(x,t)}\, d\mu(x)\,dt\\
&\quad+\int_0^\tau \int_{\G}|u(x,t)|^p\left[\rho(x)\,\partial_t \kappa(x,t)+\frac{(\kappa'(x,t))^2}{\kappa(x,t)}\right]\eta_R^2(x)\,d\mu(x)\,dt\,.
\end{aligned}
\end{equation}
Since $\gamma\ge\frac{\sigma(\sigma+1)}{\rho_0}$, we can apply Lemma \ref{lemma62} to the last term of \eqref{eq616}, it yields
$$
\int_{\G}\,\rho(x)\,|u(x,\tau)|^p{\kappa(x,\tau)}\eta_R^2(x)\,d\mu(x)\le 4\int_0^\tau\int_{\G}|u(x,t)|^p(\eta_R'(x))^2{\kappa(x,t)}\, d\mu(x)dt\,.
$$
Then, by the very definition of $\kappa$ and estimate \eqref{eta2}, we get
$$
\begin{aligned}
\int_{\G}\,\rho(x)\,|u(x,\tau)|^p&{\kappa(x,\tau)}\eta_R^2(x)\,d\mu(x)\le
4\int_0^\tau\int_{\G}|u(x,t)|^p(\eta_R'(x))^2e^{-\gamma t}[d(x,x_0)+k]^{-\sigma}\, d\mu(x)\,dt\\
&\le\frac{4C_1^2}{R^2}\int_0^\tau\int_{B_{2R}\setminus B_R}|u(x,t)|^p[d(x,x_0)+k]^{-\sigma}\,d\mu(x)\,dt.\\
&\le\frac{4C_1^2}{R^2}(R+k)^{-\sigma}\int_0^\tau\int_{B_{2R}}|u(x,t)|^p\,d\mu(x)\,dt\,.
\end{aligned}
$$
Since $\sigma=\lambda$ and due to assumption \eqref{growthcond4},
we get, for some $\tilde C>0$,
\begin{equation}\label{eq617}
\int_{\G}\,\rho(x)\,|u(x,\tau)|^p{\kappa(x,\tau)}\eta_R^2(x)\,d\mu(x)\le \frac{4C_1^2\tilde C}{R^2}\left(\frac{2R+k}{R+k}\right)^{\lambda}\le\frac{8C_1^2\tilde C}{R^2}\,.
\end{equation}
We finally take the limit as $R\to+\infty$ in \eqref{eq617}, we get
$$
\int_{\G}\,\rho(x)\,|u(x,\tau)|^p{\kappa(x,\tau)}\,d\mu(x)\le0\,.
$$
Observe that $\rho>0$ (see assumption \eqref{rhobound2}), $\kappa>0$ in $\G\times(0,T)$, therefore it necessarily holds that $u(x,\tau)=0$ for all $x\in \G$. By arbitrariness of $\tau\in(0,T]$, we get the thesis.
\medskip

We now verify item (ii). Let $1\le p<2$ and let assumption \eqref{degAss} be verified. Let $\tau\in(0,T]$, we define $v(x,t)=\eta_R(x)\,{\kappa(x,t)}$ for all $x\in \mathcal G$, $t\in(0,\tau]$.
Note that
 $$
 \Delta v(x,t)=\eta_R(x)\Delta\kappa(x,t)+2\eta_R'(x)\kappa'(x,t)+{\kappa(x,t)}\Delta \eta_R(x)\,.
 $$
We can apply Proposition \ref{prop52} with this choice of $v$, thus we get
\begin{equation}\label{eq618}
\begin{aligned}
\int_0^\tau&\int_{\G}|u(x,t)|^p\left[\rho(x)\partial_t\kappa(x,t)+\Delta \kappa(x,t)\right] \eta_R(x)\,d\mu(x)dt\\
&+\int_0^\tau\int_{\G}|u(x,t)|^p\Delta \eta_R(x) \kappa(x,t)\,d\mu(x)\,dt+2\int_0^\tau\int_{\G}|u(x,t)|^p\eta'_R(x) \kappa'(x,t)\,d\mu(x)dt\\
&\ge\int_{\G}\rho(x)|u(x,\tau)|^p\eta_R(x) \kappa(x,\tau)\,d\mu(x)+\int_0^\tau\sum_{x\in V}|u(x,t)|^p[\mathcal K(\eta_R{\kappa})](x,t)dt\,.\\
\end{aligned}
\end{equation}
Since $\gamma\ge\frac{\sigma(\sigma+1)}{\rho_0}$, we can apply Lemma \ref{lemma62} to the latter, this yields
$$
\int_0^\tau\int_{\G}|u(x,t)|^p\left[\rho(x)\partial_t\kappa(x,t)+\Delta \kappa(x,t)\right] \eta_R(x)\,d\mu(x)dt\le0
$$
Hence, \eqref{eq618} reduces to
\begin{equation}\label{eq619}
\begin{aligned}
\int_0^\tau\int_{\G}&|u(x,t)|^p\Delta \left(\eta_R(x)\right) \kappa(x,t)\,d\mu(x)dt\\
&\quad+2\int_0^\tau\int_{\G}|u(x,t)|^p\eta'_R(x)\kappa'(x,t)\,d\mu(x)dt\\
&\ge\int_{\G}\rho(x)|u(x,\tau)|^p\eta_R(x)\kappa(x,\tau)\,d\mu(x)+\int_0^\tau\sum_{x\in V}|u(x,t)|^p[\mathcal K(\eta_R\kappa)](x,t)\,dt\,.\\
&\ge\int_{\G}\rho(x)|u(x,\tau)|^p\eta_R(x)\kappa(x,\tau)\,d\mu(x)\\
&\quad+\int_0^\tau\sum_{x\in V}|u(x,t)|^p\eta_R(x)[\mathcal K({\kappa})](x,t)\,dt\\
&\quad+\int_0^\tau\sum_{x\in V}|u(x,t)|^p\kappa(x,t)[\mathcal K(\eta_R)](x)\,dt\,.\\
\end{aligned}
\end{equation}
Defining
$$
\begin{aligned}
I_{1,R}&=\int_0^\tau\int_{\G}|u(x,t)|^p\Delta \left(\eta_R(x)\right) \kappa(x,t)\,d\mu(x)dt\\
I_{2,R}&=2\int_0^\tau\int_{\G}|u(x,t)|^p\eta'_R(x) \kappa'(x,t)\,d\mu(x)dt\\
I_{3,R}&=-\int_0^\tau\sum_{x\in V}|u(x,t)|^p\kappa(x,t)[\mathcal K(\eta_R)](x)\,dt\,\\
I&=-\int_0^\tau\sum_{x\in V}|u(x,t)|^p\eta_R(x)[\mathcal K({\kappa})](x,t)\,dt\,;
\end{aligned}
$$
we can rearrange \eqref{eq619} as
\begin{equation}\label{eq620}
\int_{\G}\rho(x)|u(x,\tau)|^p\eta_R(x)\kappa(x,\tau)\,d\mu(x)\le I_{1,R}+I_{2,R}+I_{3,R}+I\,.
\end{equation}
The thesis follows if we show that $I_{i,R}\to0$ as $R\to\infty$ for $i=1,2,3$ and $I\le0$. Since $\sigma=\lambda$ and due to assumptions \eqref{growthcond4} and \eqref{eta2}, we have, for some $\tilde C>0$,
\begin{equation}\label{eq621}
\begin{aligned}
|I_{1,R}|&\le\int_0^\tau\int_{\G}|u(x,t)|^p|\Delta \eta_R(x)|e^{-\gamma t}[d(x,x_0)+k]^{-\sigma} \,d\mu(x)dt \\
&\le\frac{C_2}{R^2}\int_0^\tau\int_{B_{2R}\setminus B_R}|u(x,t)|^p [d(x,x_0)+k]^{-\sigma} \,d\mu(x)dt\\
&\le\frac{C_2}{R^2}(R+k)^{-\sigma}\int_0^\tau\int_{B_{2R}}|u(x,t)|^p \,d\mu(x)dt\le\frac{C_2}{R^2}\left(\frac{2R+k}{R+k}\right)^{\lambda} \le\frac{2 \tilde CC_2}{R^2}\,;\\
|I_{2,R}|&\le 2\lambda\int_0^\tau\int_{\G}|u(x,t)|^p|\eta_R'(x)| e^{-\gamma t}[d(x,x_0)+k]^{-\sigma-1}\,d\mu(x)dt\\
&\le \frac{2 \beta C_1}{R}\int_0^\tau\int_{B_{2R}\setminus B_R}|u(x,t)|^p[d(x,x_0)+k]^{-\sigma-1}\,d\mu(x)dt\\
&\le\frac{2\beta C_1}{R}(R+k)^{-\sigma-1}\int_0^\tau\int_{B_{2R}}|u(x,t)|^p \,d\mu(x)dt
\le \frac{4\beta \tilde C C_1}{R^2}\,.\\
\end{aligned}
\end{equation}
Moreover, from the very definition of $\mathcal K$ in \eqref{K} and $\eta_R$,
\begin{equation}\label{eq622}
\begin{aligned}
|I_{3,R}|&\le\int_0^\tau\sum_{x\in V}|u(x,t)|^p\kappa(x,t)|[\mathcal K(\eta_R)](x)|\,dt\\
&\le \frac{ C_1}{R}\int_0^\tau\sum_{x\in V}|u(x,t)|^p\kappa(x,t)\,dt\\
&\le\frac{ C_1}{R}\int_0^\tau\sum_{x\in V}|u(x,t)|^p[d(x,x_0+k]^{-\lambda}\,dt\,.
\end{aligned}
\end{equation}
Finally, we observe that, since $\kappa'\le0$, due to \eqref{e502u} and assumption \eqref{degAss}, for every $t\in(0,\tau)$,
\begin{equation*}
[\mathcal K(\kappa)](x,t)=\sum_{k=1}^{\operatorname{deg}_x^+}\kappa'(x,t)-\sum_{k=1}^{\operatorname{deg}_x^-}\kappa'(x,t)=\kappa'(x,t)\left(\operatorname{deg}_x^+-\operatorname{deg}_x^-\right)\ge0\quad\text{for every}\,\, x\in V,
\end{equation*}
therefore
\begin{equation}\label{eq623}
I=-\int_0^\tau\sum_{x\in V}|u(x,t)|^p\eta_R(x)[\mathcal K(\kappa)](x,t)dt\le0\,.
\end{equation}
By combining \eqref{eq621}, \eqref{eq622}, \eqref{eq623} and \eqref{eq620}, we get
$$
\int_{\G}\rho(x)|u(x,\tau)|^p\eta_R(x){\kappa(x,\tau)}\,d\mu(x)\le\frac{2\tilde C\,C_2}{R^2}+\frac{4 \beta\tilde C\, C_1}{kR}+\frac{ C_1}{R}\int_0^\tau\sum_{x\in V}|u(x,t)|^p[d(x,x_0+k]^{-\lambda}\,dt\,.
$$
We finally let $R\to+\infty$, we obtain
$$
\int_{\G}\rho(x)|u(x,\tau)|^p{\kappa(x,\tau)}\,d\mu(x)\le0\,.
$$
Since $\rho>0$ (see assumption \eqref{rhobound2}), $\kappa>0$ and by arbitrariness of $\tau$, the thesis follows.

\end{proof}

\bigskip
\bigskip

\noindent{\bf Acknowledgement}
The first author is funded by the Deutsche Forschungsgemeinschaft (DFG, German Research Foundation) - SFB 1283/2 2021 - 317210226.

\noindent{\bf Conflict of interest}. The authors state no conflict of interest.

\noindent{\bf Data availability statement}. There are no data associated with this research.

\end{document}